\documentclass[a4paper,11pt,reqno]{amsart}

\usepackage{latexsym,dsfont,amssymb,amsmath,amsthm}

\chardef\bslash=`\\ 

\numberwithin{equation}{section}

\newtheorem{theorem}{Theorem}[section]
\newtheorem{corollary}[theorem]{Corollary}
\newtheorem{lemma}[theorem]{Lemma}
\newtheorem{proposition}[theorem]{Proposition}

\theoremstyle{remark}
\newtheorem{remark}[theorem]{Remark}

\theoremstyle{definition}
\newtheorem{definition}[theorem]{Definition}

\newcommand\bp{\begin{proof}}
\newcommand\ep{\end{proof}}

\newcommand\3[1]{{\mathds #1}}

\newcommand{\N}{\mathbb N}
\newcommand{\Z}{\mathbb Z}
\newcommand\Zhat{\hat{\mathbb Z}}
\newcommand{\HH}{\mathbb H}
\newcommand{\Q}{\mathbb Q}
\newcommand{\R}{\mathbb R}
\newcommand{\C}{\mathbb C}
\newcommand\A{\mathbb{A}}
\newcommand\af{\mathbb{A}_f}

\newcommand{\hh}{\mathcal H}
\newcommand{\rr}{\mathcal R}

\newcommand\eps{\varepsilon}

\newcommand\bpmatrix{\begin{pmatrix}}
\newcommand\epmatrix{\end{pmatrix}}

\newcommand{\diag}[2]{\left(\begin{matrix}#1&0\\
0&#2\end{matrix}\right)}
\newcommand{\fmatr}[4]{\left(\begin{matrix}#1&#2\\
#3&#4\end{matrix}\right)}

\newcommand\mz{{\operatorname{Mat}_2}(\Z)}
\newcommand\mzp{{\operatorname{Mat}^+_2}(\Z)}

\newcommand\mr{{\operatorname{Mat}_2}(\Zhat)}
\newcommand\ma{{\operatorname{Mat}_2}(\af)}
\newcommand\mtwo{\operatorname{Mat}_2}
\newcommand\mn{{\operatorname{Mat}_n}}

\newcommand\glq{{\operatorname{GL}^+_2}(\Q)}

\newcommand\slz{{\operatorname{SL}_2}(\Z)}
\newcommand\slr{{\operatorname{SL}_2}(\Zhat)}
\newcommand\glr{{\operatorname{GL}_2}(\Zhat)}
\newcommand\gl{{\operatorname{GL}_2}}
\newcommand\glp{{\operatorname{GL}_2^+}}
\newcommand\sltwo{{\operatorname{SL}_2}}
\newcommand\gln{{\operatorname{GL}_n}}

\newcommand\G{\gl}
\newcommand\pgl{{\operatorname{PGL}^+_2}(\R)}
\newcommand\so{{\operatorname{SO}_2}(\R)}
\newcommand\pso{{\operatorname{PSO}_2}(\R)}
\newcommand\glnq{{\operatorname{GL}^+_n}(\Q)}
\newcommand\pgln{{\operatorname{PGL}^+_n}(\R)}
\newcommand\slnz{{\operatorname{SL}_n}(\Z)}

\newcommand{\hecke}[2]{\mathcal \hh({#1}, {#2})}
\newcommand{\redheck}[2]{C^*_r({#1}, {#2})}
\newcommand{\red}[1]{C^*_r(#1)}

\newcommand{\fib}[3]{#2\backslash #1\times_{#2}#3}
\newcommand{\fibb}[3]{#2\backslash #1\boxtimes_{#2}#3}

\newcommand\enu[1]{\smallskip\newline\makebox[5mm][l]{\rm(#1)}}

\begin{document}

\title{Phase transition in the Connes-Marcolli GL$_2$-system}

\thanks{Part of this work was carried out through several visits of
the first author to the Department of Mathematics at the University
of Oslo. He would like to thank the department for their
hospitality, and the SUPREMA project for the support.}

\author[M. Laca]{Marcelo Laca$^1$}
\address{Department of Mathematics and Statistics, University of
Victoria, PO Box 3045, Victoria, British Columbia, V8W 3P4, Canada.}
\email{laca@math.uvic.ca}
\thanks{$^1$) Supported by the Natural Sciences and Engineering Research
Council of Canada.}

\author[N. S. Larsen]{Nadia S. Larsen$^2$}
\address{Department of Mathematics, University of Oslo,
P.O. Box 1053 Blindern, N-0316 Oslo, Norway.}
\email{nadiasl@math.uio.no}
\thanks{$^2$) Supported by the Research Council of Norway.}

\author[S. Neshveyev]{Sergey Neshveyev$^2$}
\address{Department of Mathematics, University of Oslo,
P.O. Box 1053 Blindern, N-0316 Oslo, Norway.}
\email{sergeyn@math.uio.no}

\begin{abstract}
We develop a general framework for analyzing KMS-states on
C$^*$-algebras arising from actions of Hecke pairs. We then
specialize to the system recently introduced by Connes and Marcolli
and classify its KMS-states for inverse temperatures $\beta\ne0,1$.
In particular, we show that for each $\beta\in(1,2]$ there exists a
unique KMS$_\beta$-state.
\end{abstract}

\date{September 11, 2006; minor corrections August 27, 2007}

\maketitle

\bigskip

\section*{Introduction}

More than ten years ago Bost and Connes~\cite{bos-con} constructed
a C$^*$-dynamical system with the Galois group $G(\Q^{ab}/\Q)$ as
symmetry group and with phase transition related to properties of
zeta and $L$-functions. Since then there have been numerous, and
only partially successful, attempts to generalize the Bost-Connes
system to arbitrary number fields, see~\cite[Section~1.4]{con-mar}
for a survey. As was later emphasized by Connes, the BC-system has
yet another remarkable property: there exists a dense
$\Q$-subalgebra such that the maximal abelian extension $\Q^{ab}$
of $\Q$ arises as the set of values of a ground state of the
system on it. If one puts this property as a requirement for an
arbitrary number field, one recognizes that the problem of finding
the right analogue of the BC-system is related to Hilbert's $12$th
problem on explicit class field theory. Since the only case (in
addition to~$\Q$) for which Hilbert's problem is completely solved
is that of imaginary quadratic fields, these fields should be the
first to investigate. This has been done in recent papers of
Connes, Marcolli and
Ramachandran~\cite{con-mar,con-mar2,cmn1,cmn2}. Connes and
Marcolli~\cite{con-mar,con-mar2} constructed a $\G$-system, an
analogue of the BC-system with $\Q^*$ replaced by $\G(\Q)$. Its
specialization to a subsystem compatible with complex
multiplication in a given imaginary quadratic field gives the
right analogue of the BC-system for such a field~\cite{cmn1,cmn2}.
Later Ha and Paugam~\cite{ha-pa}, inspired by constructions of
Connes and Marcolli, proposed an analogue of the BC-system for an
arbitrary number field.

Connes and Marcolli classified KMS-states of the $\G$-system for
inverse temperatures $\beta\notin(1,2]$. It is the primary goal of
the present paper to elucidate what happens in the critical region
$(1,2]$. Along the way we develop some general tools for analyzing
systems of the type introduced by Connes and Marcolli, which can be
thought of  as crossed products of abelian algebras by Hecke
algebras.

Our approach to the problem is along the lines of that of the first
author in the case of the BC-system~\cite{lac}. Namely, in
Proposition~\ref{cm1} we show that KMS-states correspond to states
on the diagonal subalgebra which are scaled by the action of $\glq$,
or rather by the Hecke operators. As our first application we
recover in Theorem~\ref{conmar} the results of Connes and Marcolli.
We then prove our main result, Theorem~\ref{cm2}, the uniqueness of
a KMS$_\beta$-state for each $\beta\in(1,2]$. The strategy is
similar to that of the third author in the BC-case~\cite{nes}.
Namely, we prove the uniqueness and ergodicity, under the action of
$\glq$, of the measure defining a symmetric KMS$_\beta$-state by
analyzing an explicit formula for the projection onto the space of
$\mzp$-invariant functions, see Lemma~\ref{symm} and
Corollary~\ref{glqergodicpglma}, and then derive from this the main
uniqueness result. There are two main complications compared to the
BC-case. The first is that instead of semigroup actions we now have
to deal with representations of Hecke algebras. The second is the
presence in the system of a continuous component corresponding to
the infinite place. As a result, the critical step now is to prove
the uniqueness of a symmetric, that is, $\G(\Zhat)$-invariant,
KMS$_\beta$-state, while in the BC-case the analogous statement is
almost obvious. To show this uniqueness we use a deep result of
Clozel, Oh and Ullmo~\cite{cou} on equidistribution of Hecke points.
We point out that, as opposed to the BC-case, there are many
symmetric states for $\beta>2$, which can be easily seen from
Theorem~\ref{conmar} below.

\bigskip

\section{Proper actions and groupoid C$^*$-algebras}
\label{sproper}

Let $G$ be a countable group acting on a locally compact second
countable space $X$. The reduced crossed product $C_0(X)\rtimes_rG$
is the reduced C$^*$-algebra of the transformation groupoid $G\times
X$ with unit space $X$, source and range maps $(g,x)\mapsto x$ and
$(g,x)\mapsto gx$, respectively, and the product
$$
(g,hx)(h,x)=(gh,x).
$$
If the restriction of the action to a subgroup $\Gamma$ of $G$ is
free and proper, we can introduce a new groupoid
$\fib{G}{\Gamma}{X}$ by taking the quotient of $G\times X$ by the
action of $\Gamma\times\Gamma$ defined by
\begin{equation} \label{eaction}
(\gamma_1,\gamma_2)(g,x)=(\gamma_1g\gamma_2^{-1},\gamma_2x).
\end{equation}
Thus the unit space of $\fib{G}{\Gamma}{X}$ is $\Gamma\backslash X$,
and the product is induced from that on $G\times X$. This groupoid
is Morita equivalent in the sense of~\cite{mrj} to the
transformation groupoid $G\times X$. Although we will not need this
result, let us briefly recall the argument. By definition of Morita
equivalence first of all we have to find a space $Z$ with commuting
actions of our groupoids. We take $Z=G\times_\Gamma X$, the quotient
of $G\times X$ by the action of $\Gamma$ given by
$\gamma(g,x)=(g\gamma^{-1}, \gamma x)$. The left and right actions
of the groupoid $G\times X$ on itself induce a left action of
$G\times X$ and a right action of $\fib{G}{\Gamma}{X}$ on~$Z$. The
map $Z\to\Gamma\backslash X$, $\Gamma(g,x)\mapsto\Gamma x$, induces
a homeomorphism between the quotient of $Z$ by the action of
$G\times X$ and the unit space $\Gamma\backslash X$ of the groupoid
$\fib{G}{\Gamma}{X}$. Similarly, the map $Z\to X$,
$\Gamma(g,x)\mapsto gx$, induces a homeomorphism between the
quotient of $Z$ by $\fib{G}{\Gamma}{X}$ and $X$. Thus the groupoids
are indeed Morita equivariant. Recall then that
by~\cite[Theorem~2.8]{mrj} the corresponding reduced C$^*$-algebras
are Morita equivalent.

If the action of $\Gamma$ is proper but not free, the quotient space
$\fib{G}{\Gamma}{X}$ is no longer a groupoid, since the composition
of classes using representatives will in general depend on the
choice of representatives. As was observed in~\cite{con-mos} and
\cite{con-mar}, nevertheless, the same formula for convolution of
two functions as in the groupoid case gives us a well-defined
algebra, and by completion we get a C$^*$-algebra. In more detail,
consider the space $C_c(\fib{G}{\Gamma}{X})$ of continuous compactly
supported functions on $\fib{G}{\Gamma}{X}$. We consider its
elements as $(\Gamma\times\Gamma)$-invariant functions on $G\times
X$, and define a convolution of two such functions by
\begin{equation} \label{econv}
(f_1*f_2)(g,x)=\sum_{h\in\Gamma\backslash G}f_1(gh^{-1},hx)
f_2(h,x).
\end{equation}
To see that the convolution is well-defined, assume the support of
$f_i$ is contained in $(\Gamma\times\Gamma)(\{g_i\}\times U_i)$,
where $g_i\in G$ and $U_i$ is a compact subset of $X$. Let
$\{\gamma_1,\dots,\gamma_n\}$ be the set of all elements
$\gamma\in\Gamma$ such that $\gamma g_2U_2\cap U_1\ne\emptyset$.
Note that this set is finite since the action of $\Gamma$ is assumed
to be proper. If $f_2(h,x)\ne0$ then there exists $\gamma\in\Gamma$
such that $h\gamma^{-1}\in\Gamma g_2$ and $\gamma x\in U_2$. Since
the number of $\gamma$'s such that $\gamma x\in U_2$ is finite, we
already see that the sum above is finite. If furthermore
$f_1(gh^{-1},hx)\ne0$ then replacing $h$ by another representative
of the right coset $\Gamma h$ we may assume that $gh^{-1}\in\Gamma
g_1$ and $hx\in U_1$. Then if $h\gamma^{-1}=\tilde\gamma g_2$ with
$\tilde\gamma\in\Gamma$, we get $hx=\tilde\gamma g_2\gamma
x\in\tilde\gamma g_2U_2$. Hence $\tilde\gamma=\gamma_i$ for
some~$i$, and therefore $g\in\Gamma g_1h=\Gamma g_1\gamma_i
g_2\gamma$. Thus the support of $f_1*f_2$ is contained in the union
of the sets $(\Gamma\times\Gamma)(\{g_1\gamma_ig_2\}\times U_2)$, so
$f_1*f_2\in C_c(\fib{G}{\Gamma}{X})$ and the latter space becomes an
algebra. It is not difficult to check that the convolution is
associative.

Define also an involution on $C_c(\fib{G}{\Gamma}{X})$ by
\begin{equation} \label{einv}
f^*(g,x)=\overline{f((g,x)^{-1})}
=\overline{f(g^{-1},gx)}.
\end{equation}
If the support of $f$ is contained in
$(\Gamma\times\Gamma)(\{g_0\}\times U)$ for $g_0\in G$ and compact
$U\subset X$, then the support of $f^*$ is contained in
$$
((\Gamma\times\Gamma)(\{g_0\}\times U))^{-1}=
(\Gamma\times\Gamma)(\{g_0\}\times U)^{-1}=
(\Gamma\times\Gamma)(\{g_0^{-1}\}\times g_0U),
$$
so indeed $f^*\in C_c(\fib{G}{\Gamma}{X})$.

For each $x\in X$ we define a $*$-representation $\pi_x\colon
C_c(\fib{G}{\Gamma}{X})\to B(\ell^2(\Gamma\backslash G))$ by
\begin{equation} \label{erep}
\pi_x(f)\delta_{\Gamma h}=\sum_{g\in\Gamma\backslash
G}f(gh^{-1},hx)\delta_{\Gamma g},
\end{equation}
where $\delta_{\Gamma g}$ denotes the characteristic function of
the coset $\Gamma g$. It is standard to show that the operators
$\pi_x(f)$ are bounded, but we include a proof for the reader's
convenience.

\begin{lemma}
For each $f\in C_c(\fib{G}{\Gamma}{X})$ the operators $\pi_x(f)$,
$x\in X$, are uniformly bounded.
\end{lemma}

\bp For $\xi_1,\xi_2\in\ell^2(\Gamma\backslash G)$ we have
\begin{eqnarray*}
|(\pi_x(f)\xi_1,\xi_2)|&\le&\sum_{g,h\in\Gamma\backslash G}
|f(gh^{-1},hx)|\,|\xi_1(h)|\,|\xi_2(g)|\\
&\le&\left(\sum_{g,h\in\Gamma\backslash G}
|f(gh^{-1},hx)|\,|\xi_1(h)|^2\right)^{1/2}
\left(\sum_{g,h\in\Gamma\backslash G}
|f(gh^{-1},hx)|\,|\xi_2(g)|^2\right)^{1/2}.
\end{eqnarray*}
Thus if we denote by $\|f\|_I$ the quantity
$$
\max\left\{\sup_{x\in X,\,h\in G}\sum_{g\in\Gamma\backslash G}
|f(gh^{-1},hx)|,\ \sup_{x\in X,\,g\in G}\sum_{h\in\Gamma\backslash
G} |f(gh^{-1},hx)|\right\},
$$
we get $\|\pi_x(f)\|\le\|f\|_I$ for any $x\in X$, so it suffices to
show that $\|f\|_I$ is finite. Replacing $x$ by $h^{-1}x$ and $g$
by $gh$ in the first supremum above, we see that this supremum
equals
$$
\|f\|_{I,s}:=\sup_{x\in X}\sum_{g\in\Gamma\backslash G} |f(g,x)|.
$$
Observe next that $f(gh^{-1},hx)=\overline{f^*(hg^{-1},gx)}$, so
that the second supremum is equal to $\|f^*\|_{I,s}$. Therefore
$\|f\|_I=\max\{\|f\|_{I,s},\|f^*\|_{I,s}\}$. It remains to show
that $\|f\|_{I,s}$ is finite for any $f\in
C_c(\fib{G}{\Gamma}{X})$.

Assume the support of $f$ is contained in
$(\Gamma\times\Gamma)(\{g_0\}\times U)$ for some $g_0\in G$ and
compact $U\subset X$. Since the action of $\Gamma$ is proper,
there exists $n\in\N$ such that the sets $\gamma_i U$,
$i=1,\dots,n+1$, have trivial intersection for any different
$\gamma_1,\dots,\gamma_{n+1}\in\Gamma$. Now if $f(g,x)\ne0$ for
some $g$ and $x$, there exists $\gamma\in\Gamma$ such that
$g\gamma^{-1}\in\Gamma g_0$ and $\gamma x\in U$. Since the number
of $\gamma$'s such that $\gamma x\in U$ is at most~$n$, we see
that for each $x\in X$ the sum in the definition of $\|f\|_{I,s}$
has at most $n$ nonzero summands. Hence $\|f\|_{I,s}$ is finite,
and the proof of the lemma is complete. \ep

We denote by $\red{\fib{G}{\Gamma}{X}}$ the completion of
$C_c(\fib{G}{\Gamma}{X})$ in the norm defined by the
representation $\oplus_{x\in X}\pi_x$, that is,
$$
\|f\|=\sup_{x\in X}\|\pi_x(f)\|.
$$
Denoting by $U_g$ the unitary operator on $\ell^2(\Gamma\backslash
G)$ such that $U_g\delta_{\Gamma h}=\delta_{\Gamma hg^{-1}}$, we get
$U_g\pi_x(f)U^*_g=\pi_{gx}(f)$. Hence $\|\pi_x(f)\|=\|\pi_{gx}(f)\|$
and so the supremum above is actually over $G\backslash X$.

\smallskip

Using the embedding $X\hookrightarrow G\times X$, $x\mapsto
(e,x)$, we may consider $\Gamma\backslash X$ as an open subset of
$\fib{G}{\Gamma}{X}$, and then the algebra $C_0(\Gamma\backslash
X)$ as a subalgebra of $\red{\fib{G}{\Gamma}{X}}$. More generally,
any bounded continuous function on $\Gamma\backslash X$ defines
a multiplier of $\red{\fib{G}{\Gamma}{X}}$.

\begin{lemma} \label{Cond}
There exists a conditional expectation $E\colon
\red{\fib{G}{\Gamma}{X}} \to C_0(\Gamma\backslash X)$ such that
$$
E(f)(x)=f(e,x)\ \ \text{for}\ \ f\in C_c(\fib{G}{\Gamma}{X}).
$$
\end{lemma}

\bp For each $x\in X$ define a state $\omega_x$ on
$\red{\fib{G}{\Gamma}{X}}$ by
$$
\omega_x(a)=(\pi_x(a)\delta_\Gamma,\delta_\Gamma).
$$
Then the function $E(a)$ on $X$ defined by $E(a)(x)=\omega_x(a)$ is
bounded by $\|a\|$. Since $E(f)(x)=f(e,x)$ for $f\in
C_c(\fib{G}{\Gamma}{X})$, we conclude that $E(a)\in
C_0(\Gamma\backslash X)$ for every $a\in \red{\fib{G}{\Gamma}{X}}$.
Thus $E$ is the required conditional expectation. \ep

Let $Y\subset X$ be a $\Gamma$-invariant clopen subset. Then, as
we already observed, the characteristic function
$\31_{\Gamma\backslash Y}$ of the set $\Gamma\backslash Y$ is an
element of the multiplier algebra of $\red{\fib{G}{\Gamma}{X}}$.
Denote by $\fibb{G}{\Gamma}{Y}$ the quotient of the space
$$
\{(g,x)\mid g\in G,\,x\in Y,\, gx\in Y\}
$$
by the action of $\Gamma\times\Gamma$ defined as in (\ref{eaction}).
Then
$$
\31_{\Gamma\backslash
Y}C_c(\fib{G}{\Gamma}{X})\31_{\Gamma\backslash Y}
=C_c(\fibb{G}{\Gamma}{Y}).
$$
Therefore the algebra $\31_{\Gamma\backslash
Y}\red{\fib{G}{\Gamma}{X}}\31_{\Gamma\backslash Y}$, which we
shall denote by $\red{\fibb{G}{\Gamma}{Y}}$, is a completion of
the algebra of compactly supported functions on
$\fibb{G}{\Gamma}{Y}$ with convolution product given by
$$
(f_1*f_2)(g,y)=\sum_{h\in\Gamma\backslash G\colon hy\in
Y}f_1(gh^{-1},hy) f_2(h,y),
$$
and involution
$$
f^*(g,y)=\overline{f(g^{-1},gy)}.
$$
Note that $\pi_x(\31_{\Gamma\backslash Y})$ is the projection onto
the subspace $\ell^2(\Gamma\backslash G_x)$ of
$\ell^2(\Gamma\backslash G)$, where the subset $G_x$ of $G$ is
defined by
$$
G_x=\{g\in G\,\mid\, gx\in Y\},
$$
and then
$$
\pi_x(f)\delta_{\Gamma h}=\sum_{g\in\Gamma\backslash
G_x}f(gh^{-1},hx)\delta_{\Gamma g}
$$
for $h\in G_x$ and $f\in C_c(\fibb{G}{\Gamma}{Y})$. In particular,
$\pi_x(f)=0$ if $x\notin GY$. As we already remarked, the
representations $\pi_x$ and $\pi_{gx}$ are unitarily equivalent for
any $g\in G$. Thus we may conclude that $\red{\fibb{G}{\Gamma}{Y}}$
is precisely the completion of $C_c(\fibb{G}{\Gamma}{Y})$ in the
norm
$$
\|f\|=\sup_{y\in Y}\|\pi_y(f)\|.
$$
This is how the algebra $\red{\fibb{G}{\Gamma}{Y}}$ was defined (in
a particular case) in~\cite[Proposition~1.23]{con-mar}.

\smallskip

Returning to the algebra $\red{\fib{G}{\Gamma}{X}}$, our next goal
is to show that under an extra assumption its multiplier algebra
contains other interesting elements in addition to the
$\Gamma$-invariant functions on~$X$.

Recall that $(G,\Gamma)$ is called a Hecke pair if $\Gamma$ and
$g\Gamma g^{-1}$ are commensurable for any $g\in G$, that is,
$\Gamma\cap g\Gamma g^{-1}$ is a subgroup of $\Gamma$ of finite
index. Equivalently, every double coset of $\Gamma$ contains
finitely many right (and left) cosets of $\Gamma$, so that
$$
R_\Gamma(g):=|\Gamma\backslash\Gamma g\Gamma|<\infty\ \text{for
any}\ g\in G.
$$
Then the space $\hecke{G}{\Gamma}$ of finitely supported functions
on $\Gamma\backslash G/\Gamma$ is a $*$-algebra with  product
$$
(f_1*f_2)(g)=\sum_{h\in\Gamma\backslash G}f_1(gh^{-1})f_2(h)
$$
and involution $f^*(g)=\overline{f(g^{-1})}$, see e.g.~\cite{kri}.
This algebra is represented on $\ell^2(\Gamma\backslash G)$ by
$$
f\delta_{\Gamma h}=\sum_{g\in\Gamma\backslash
G}f(gh^{-1})\delta_{\Gamma g},
$$
see \cite{bos-con}. The corresponding completion is called the
reduced Hecke C$^*$-algebra of $(G,\Gamma)$ and denoted by
$\redheck{G}{\Gamma}$. We shall denote by $[g]$ the characteristic
function of the double coset $\Gamma g\Gamma$ considered as an
element of the Hecke algebra.

We may consider elements of $\hecke{G}{\Gamma}$ as continuous
functions on $\fib{G}{\Gamma}{X}$. Although these functions are not
compactly supported in general, the formulas defining the
$*$-algebra structure and the regular representation of
$\hecke{G}{\Gamma}$ coincide with (\ref{econv})-(\ref{erep}).
Furthermore, the convolution of an element of $\hecke{G}{\Gamma}$
with a compactly supported function on $\fib{G}{\Gamma}{X}$ gives a
compactly supported function. Indeed, if $f_1=[g_1]$ and the support
of $f_2\in C_c(\fib{G}{\Gamma}{X})$ is contained in
$(\Gamma\times\Gamma)(\{g_2\}\times U)$ for a compact $U\subset X$,
then the support of $f_1*f_2$ is contained in
$(\Gamma\times\Gamma)(g_1\Gamma g_2\times U)$. Since
$\Gamma\backslash \Gamma g_1\Gamma g_2$ is finite, we see that
$f_1*f_2$ is compactly supported on $\fib{G}{\Gamma}{X}$. We may
therefore conclude the following.

\begin{lemma} \label{HeckeM}
If $(G,\Gamma)$ is a Hecke pair, then the reduced Hecke
C$^*$-algebra $\redheck{G}{\Gamma}$ is contained in the multiplier
algebra of the C$^*$-algebra $\red{\fib{G}{\Gamma}{X}}$.
\end{lemma}

It is then tempting to think of $\red{\fib{G}{\Gamma}{X}}$ as a
crossed product of $C_0(\Gamma\backslash X)$ by an action of the
Hecke pair $(G,\Gamma)$. This point of view has been formalized by
Tzanev~\cite{tza} who introduced a notion of a crossed product of
an algebra by an action of a Hecke pair.

\begin{remark} \label{heckecross}
We defined $\red{\fib{G}{\Gamma}{X}}$ assuming that the action of
$\Gamma$ on $X$ is proper. It is however easy to see that the
construction makes sense under the following weaker assumptions:
$\fib{G}{\Gamma}{X}$ is Hausdorff, and if for a compact set
$K\subset X$ we put $\Gamma_K=\{\gamma\in\Gamma\mid\gamma K\cap
K\ne\emptyset\}$ then the set $\Gamma\backslash\Gamma g\Gamma_K$
is finite for any $g\in G$. Note that the second assumption is
automatically satisfied when $(G,\Gamma)$ is a Hecke pair.
\end{remark}

\bigskip

\section{Dynamics and KMS-states}

Assume as above that we have an action of $G$ on $X$ such that the
action of $\Gamma\subset G$ is proper, and $Y\subset X$ is a
$\Gamma$-invariant clopen set. Assume now that we are given a
homomorphism
$$
N\colon G\to\R^*_+=(0,+\infty)
$$
such that $\Gamma$ is contained in the kernel of $N$. Then we
define a one-parameter group of automorphisms of $\red{\fib G
\Gamma X}$ by
$$
\sigma_t(f)(g,x)=N(g)^{it}f(g,x)\ \ \text{for}\ \ f\in C_c(\fib G
\Gamma X).
$$
More precisely, if we denote by $\bar N$ the selfadjoint operator
on $\ell^2(\Gamma\backslash G)$ defined by
$$
\bar N\delta_{\Gamma g}=N(g)\delta_{\Gamma g},
$$
then the dynamics $\sigma_t$ is spatially implemented by the
unitary operator $\oplus_{x\in X}\bar N^{it}$ on $\oplus_{x\in
X}\ell^2(\Gamma\backslash G)$. In other words,
$$
\pi_x(\sigma_t(a))=\bar N^{it}\pi_x(a)\bar N^{-it}\ \ \text{for
all}\ \ x\in X.
$$

Recall, see e.g.~\cite{kus}, that a semifinite $\sigma$-invariant
weight $\varphi$ is called a $\sigma$-KMS$_\beta$-weight if
$$
\varphi(aa^*)=\varphi(\sigma_{i\beta/2}(a)^*\sigma_{i\beta/2}(a))
$$
for any $\sigma$-analytic element $a$. The following result will
be the basis of our analysis of KMS-weights.

\begin{proposition} \label{Scaling1}
Assume the action of $G$ on $X$ is free, so that in particular
$\fibb{G}{\Gamma}{Y}$ is a genuine groupoid. Then for any
$\beta\in\R$ there exists a one-to-one correspondence between
$\sigma$-KMS$_\beta$ weights $\varphi$ on $\red{\fibb G \Gamma Y}$
with domain of definition containing $C_c(\Gamma\backslash Y)$ and
Radon measures $\mu$ on $Y$ such that
\begin{equation} \label{escaling}
\mu(gZ)=N(g)^{-\beta}\mu(Z)
\end{equation}
for every $g\in G$ and every compact subset $Z\subset Y$ such that
$gZ\subset Y$. Namely, such a measure $\mu$ is $\Gamma$-invariant,
so it determines a measure $\nu$ on $\Gamma\backslash Y$ such that
\begin{equation} \label{equotientm}
\int_Y f(y)d\mu(y)=\int_{\Gamma\backslash Y}\left(\sum_{y\in
p^{-1}(\{t\})}f(y)\right)d\nu(t)\ \ \text{for}\ \ f\in C_c(Y),
\end{equation}
where $p\colon Y\to\Gamma\backslash Y$ is the quotient map, and the
associated weight $\varphi$ is given by
$$
\varphi(a)=\int_{\Gamma\backslash Y}E(a)(x)d\nu(x),
$$
where $E$ is the conditional expectation defined
in Lemma~\ref{Cond}.
\end{proposition}

\bp For $\Gamma=\{e\}$ the result is well-known, see
e.g.~\cite[Proposition~II.5.4]{ren}. For arbitrary $\Gamma$ the
result can be deduced from the fact that the C$^*$-algebra
$\red{\fibb G \Gamma Y}$ is Morita equivalent to the C$^*$-algebra
$\31_{Y}(C_0(X)\rtimes_r G)\31_{Y}$ and general results on
KMS-weights on Morita equivalent algebras,
see~\cite[Theorem~3.2]{LN}. However, a more elementary way is to
argue as follows.

Since the action of $\Gamma$ on $Y$ is free, the quotient space
$\fibb G \Gamma Y$ is an etale groupoid. In fact it is an etale
equivalence relation on $\Gamma\backslash Y$, or an $r$-discrete
principal groupoid in the terminology of~\cite{ren}. To see this we
have to check that the isotropy group of every point in
$\Gamma\backslash Y$ is trivial, that is, if $g\in G$ is such that
$gy\in Y$ and $p(gy)=p(y)$ for some $y\in Y$ then $(g,y)$ belongs to
the $(\Gamma\times\Gamma)$-orbit of~$(e,y)$. But if $p(gy)=p(y)$,
there exists $\gamma\in\Gamma$ such that $\gamma gy=y$. Then $\gamma
g=e$, since the action of~$G$ is free, and therefore
$(g,y)=(\gamma^{-1},e)(e,y)$.

It is then standard to show using \cite[Proposition~II.5.4]{ren}
that $\sigma$-KMS$_\beta$ weights (with domain of definition
containing $C_c(\Gamma\backslash Y)$) on the C$^*$-algebra
$\red{\fibb G \Gamma Y}$ of the etale equivalence relation are in
one-to-one correspondence with measures $\nu$ on $\Gamma\backslash
Y$ with Radon-Nikodym cocycle $(p(y),p(gy))\mapsto N(g)^\beta$. The
latter means the following, see~\cite[Definition~I.3.4]{ren}. Assume
$Y_0$ is an open subset of~$Y$ such that the map $p\colon
Y\to\Gamma\backslash Y$ is injective on $Y_0$, and $g\in G$ is such
that $gY_0\subset Y$. Define an injective map $\tilde g\colon
p(Y_0)\to p(gY_0)$ by $\tilde gp(y)=p(gy)$ for $y\in Y_0$, and let
$\tilde g_*\nu$ be the push-forward of the measure $\nu$ under the
map $\tilde g$, that is, $\tilde g_*\nu(Z)=\nu(\tilde g^{-1}(Z))$
for $Z\subset p(gY_0)$. Then
$$
\frac{d\tilde g_*\nu}{d\nu}=N(g)^{\beta}\ \text{on}\ p(gY_0).
$$
If we denote by $\mu$ the $\Gamma$-invariant measure on $Y$
corresponding to $\nu$ via (\ref{equotientm}), then to say that the
Radon-Nikodym cocycle of $\nu$ is $(p(y),p(gy))\mapsto N(g)^\beta$
is the same as saying that $\mu$ satisfies the scaling condition
(\ref{escaling}). \ep

It will be convenient to extend the measure $\mu$ to the set $GY$.

\begin{lemma} \label{extensionX}
If $\mu$ is a measure on $Y$ as in Proposition~\ref{Scaling1}, then
it extends uniquely to a Radon measure on $GY\subset X$ satisfying
(\ref{escaling}) for $Z\subset GY$ and $g\in G$.
\end{lemma}

\bp A more general result on extensions of KMS-weights is proved
in~\cite{LN}, but the present particular case has the following
elementary proof. Choose Borel subsets $Y_i\subset Y$ and elements
$g_i\in G$ such that $GY$ is the disjoint union of the sets
$g_i^{-1}Y_i$. There is only one choice for a measure extending
$\mu$ and satisfying~(\ref{escaling}) on $GY$, namely, for a Borel
subset $Z\subset GY$ let
$$
\mu(Z)=\sum_iN(g_i)^\beta\mu(g_iZ\cap Y_i).
$$
To show that $\mu(Z)$ is independent of any choices and that the
extension satisfies (\ref{escaling}), assume $GY$ is a disjoint
union of sets $h_j^{-1}Z_j$ for some $h_j\in G$ and Borel
$Z_j\subset Y$. Let $g\in G$. Then
\begin{eqnarray*}
\sum_iN(g_i)^\beta\mu(g_igZ\cap Y_i)
&=&\sum_iN(g_i)^\beta\sum_j\mu(g_igZ\cap Y_i\cap g_igh_j^{-1}Z_j)\\
&=&\sum_iN(g_i)^\beta\sum_jN(g_igh^{-1}_j)^{-\beta}\mu(h_jZ\cap
h_jg^{-1}g^{-1}_iY_i\cap Z_j)\\
&=&N(g)^{-\beta}\sum_jN(h_j)^\beta
\sum_i\mu(h_jZ\cap h_jg^{-1}g^{-1}_iY_i\cap Z_j)\\
&=&N(g)^{-\beta}\sum_jN(h_j)^\beta\mu(h_jZ\cap Z_j).
\end{eqnarray*}
Taking $g=e$ we see that the extension of $\mu$ to $GY$ is
well-defined. But then for arbitrary $g$ the above identity reads
as $\mu(gZ)=N(g)^{-\beta}\mu(Z)$. \ep

\begin{remark} \label{vonneumann}
In the notation of Proposition~\ref{Scaling1} choose a
$\mu$-measurable subset $U$ of $Y$ such that $p\colon
Y\to\Gamma\backslash Y$ is injective on $U$ and
$p(U)=\Gamma\backslash Y$. Then the map $p$ induces an isomorphism
between the restriction $\rr_{G,U}$ of the $G$-orbit equivalence
relation on $X$ to $U$ and the principal groupoid
$\fibb{G}{\Gamma}{Y}$. Hence
$\pi_\varphi(\red{\fibb{G}{\Gamma}{Y}})''$ is isomorphic to the
von Neumann algebra $W^*(\rr_{G,U},\mu)$ of~$(\rr_{G,U},\mu)$,
see~\cite{fel-moo}. Extend the measure $\mu$ to a
$G$-quasi-invariant measure on $GY$, which we still denote by
$\mu$. Then $W^*(\rr_{G,U},\mu)$ is the reduction of the von
Neumann algebra of the $G$-orbit equivalence relation on
$(GY,\mu)$ by the projection~$\31_U$. Therefore
$$
\pi_\varphi(\red{\fibb{G}{\Gamma}{Y}})''\cong
\31_U(L^\infty(GY,\mu)\rtimes G)\31_U.
$$
\end{remark}

In some cases an argument similar to the proof of
Lemma~\ref{extensionX} allows us to describe all measures
satisfying~(\ref{escaling}).

\begin{lemma} \label{extension}
Let $Y_0$ be a $\Gamma$-invariant Borel subset of $Y$ such that
\enu{i} if $gY_0\cap Y_0\ne\emptyset$ for some $g\in G$ then
$g\in\Gamma$; \enu{ii} for any $y\in Y$ there exists $g\in G$ such
that $gy\in Y_0$.

Then any $\Gamma$-invariant Borel measure on $Y_0$ extends
uniquely to a Borel measure on $Y$ satisfying~(\ref{escaling}).
\end{lemma}

\bp Let $\mu_0$ be a $\Gamma$-invariant measure on $Y_0$. Since the
assumptions imply that $Y$ is a disjoint union of translates of
$Y_0$ by representatives of the right cosets of $\Gamma$, that is,
$Y=\sqcup_{h\in\Gamma\backslash G}(h^{-1}Y_0\cap Y)$, there is only
one choice for a measure $\mu$ extending $\mu_0$ and
satisfying~(\ref{escaling}), namely,
$$
\mu(Z)=\sum_{h\in\Gamma\backslash G}N(h)^\beta\mu_0(hZ\cap Y_0).
$$
Since $\mu_0$ is $\Gamma$-invariant, $\mu(Z)$ is independent of
the choice of representatives, so all we need to check is
that~(\ref{escaling}) holds. Let $g\in G$. Then
$$
\mu(gZ)=\sum_{h\in\Gamma\backslash G}N(h)^\beta\mu_0(hgZ\cap Y_0)
=N(g)^{-\beta}
\sum_{h\in\Gamma\backslash G}N(hg)^\beta\mu_0(hgZ\cap Y_0)
=N(g)^{-\beta}\mu(Z),
$$
and the proof is complete. \ep

Although the condition for a measure $\nu$ on $\Gamma\backslash Y$
to define a KMS-weight is easier to formulate in terms of the
corresponding $\Gamma$-invariant measure on $Y$, it will also be
important to work directly with $\nu$. For this we introduce the
following operators on functions on $\Gamma\backslash X$. We shall
often consider functions on $\Gamma\backslash X$ as
$\Gamma$-invariant functions on $X$.

\begin{definition}
Let $G$ act on a set $X$ and suppose $(G,\Gamma)$ is a Hecke a
pair. The {\it Hecke operator} associated to $g\in G$ is the
operator $T_g$ on $\Gamma$-invariant functions on $X$ defined by
$$
(T_gf)(x)=\frac{1}{R_\Gamma(g)}\sum_{h\in\Gamma\backslash\Gamma
g\Gamma}f(hx).
$$
\end{definition}

Clearly $T_gf$ is again $\Gamma$-invariant. It is not difficult to
check that the map $[g^{-1}]\to R_\Gamma(g)T_g$ is a
representation of the Hecke algebra $\hecke G \Gamma$ on the space
of $\Gamma$-invariant functions (notice that for $X=G$ this is
exactly the way we defined the regular representation of
$\hecke{G}{\Gamma}$, so by decomposing an arbitrary $X$ into
$G$-orbits one can obtain the general case without any
computations).

\smallskip

The following three lemmas will be our main computational tools.

\begin{lemma} \label{scaling2}
Suppose $\mu$ is as in Proposition~\ref{Scaling1} and $\nu$ is the
measure on $\Gamma\backslash Y$ determined by (\ref{equotientm}).
Assume further that $Y=X$, the action of $G$ on $X$ is free and that
$(G,\Gamma)$ is a Hecke pair with modular function
$\Delta_\Gamma(g):=R_\Gamma(g^{-1})/R_\Gamma(g)$. Then for any
positive measurable function $f$ on $\Gamma\backslash X$ and $g\in
G$ we have
$$
\int_{\Gamma\backslash X}T_gfd\nu=\Delta_\Gamma(g)N(g)^\beta
\int_{\Gamma\backslash X}fd\nu.
$$
\end{lemma}

\bp Fix a point $x\in X$. We claim that there exists a
neighbourhood $U$ of $x$ such that the sets $hU$ are disjoint for
different $h$ in $\Gamma g^{-1}\Gamma$. Indeed, choose
representatives $h_1,\dots,h_n$ of the right $\Gamma$-cosets
contained in $\Gamma g^{-1}\Gamma$. Since the action of $\Gamma$
is proper, there exists a neighbourhood~$U$ of $x$ such that if
$h_iU\cap\gamma h_jU\ne\emptyset$ for some $i$, $j$ and
$\gamma\in\Gamma$ then $h_ix=\gamma h_jx$. But since the action of
$G$ is free, the latter equality is possible only when $h_i=\gamma
h_j$, so that $i=j$ and $\gamma=e$. Thus $h_iU\cap\gamma
h_jU=\emptyset$ if $i\ne j$ or $\gamma\ne e$. Since $\Gamma
g^{-1}\Gamma=\cup^n_{k=1}\Gamma h_k$, this proves the claim.

The set $\Gamma g^{-1}\Gamma U$ is therefore a disjoint union of
the sets $hU$, $h\in \Gamma g^{-1}\Gamma$. So we can write
$$
\sum_{h\in\Gamma\backslash\Gamma g\Gamma}\31_{h^{-1}\Gamma U}
=\31_{\Gamma g^{-1}\Gamma U}=\sum_{h\in\Gamma\backslash\Gamma
g^{-1}\Gamma}\31_{\Gamma hU},
$$
Denoting by $p\colon X\mapsto\Gamma\backslash X$ the quotient map,
we can rewrite the above in terms of functions on
$\Gamma\backslash X$ as
$$
R_\Gamma(g)T_g(\31_{p(U)})=\31_{p(\Gamma g^{-1}\Gamma U)}
=\sum_{h\in\Gamma\backslash\Gamma g^{-1}\Gamma}\31_{p(hU)}.
$$
It follows that
$$
R_\Gamma(g)\int_{\Gamma\backslash X}T_g(\31_{p(U)})d\nu
=\sum_{h\in\Gamma\backslash\Gamma g^{-1}\Gamma}\nu(p(hU))
=\sum_{h\in\Gamma\backslash\Gamma g^{-1}\Gamma}\mu(hU)=
R_\Gamma(g^{-1})N(g)^\beta\nu(p(U)).
$$
In other words, the identity in the lemma holds for
$f=\31_{p(U)}$. Since this is true for any $x$ and sufficiently
small neighbourhood $U$ of $x$, we get the result. \ep

Notice that by applying the above lemma to the characteristic
function of $X$ we get the following: if a group $G$ acts freely on
a space $X$ with a $G$-invariant measure $\mu$, and $\Gamma$ is an
almost normal subgroup of $G$ (that is, $(G,\Gamma)$ is a Hecke
pair) such that the action of $\Gamma$ on $X$ is proper and
$0<\mu(\Gamma\backslash X)<\infty$, then $\Delta_\Gamma(g)=1$ for
any $g\in G$. The same is true if we assume that the action of~$G$
on $(X,\mu)$ is only essentially free.

\begin{lemma} \label{measure}
Suppose $\mu$ is as in Proposition~\ref{Scaling1} and $\nu$ is the
measure on $\Gamma\backslash Y$ determined by (\ref{equotientm}).
Assume the action of $G$ on~$X$ is free and that $(G,\Gamma)$ is a
Hecke pair. Assume further that $Y_0$ is a $\Gamma$-invariant
measurable subset of~$Y$ such that if $gY_0\cap Y_0\ne\emptyset$
for some $g\in G$ then $g\in\Gamma$. Then for any $g\in G$ such
that $gY_0\subset Y$, measurable $Z\subset\Gamma\backslash Y_0$
and positive measurable function $f$ on $\Gamma\backslash Y$ we
have
$$
\int_{\Gamma g Z}fd\nu=N(g)^{-\beta} R_\Gamma(g)\int_ZT_gfd\nu,
$$
where $\Gamma gZ=p(\Gamma gp^{-1}(Z))$ and $p\colon
X\to\Gamma\backslash X$ is the quotient map. In particular,
$\nu(\Gamma gZ)=N(g)^{-\beta}R_\Gamma(g)\nu(Z)$.
\end{lemma}

\bp Suppose $Z\subset\Gamma\backslash Y_0$ is measurable, and
choose $U\subset Y_0$ measurable such that $Z=p(U)$ and $p$ is
injective on $U$. For $g\in G$ let $h_1,\dots,h_n$ be
representatives of the right $\Gamma$-cosets contained in $\Gamma
g\Gamma$. We claim that the map $p$ is injective on
$h_1U,\dots,h_nU$, and the images of these sets are disjoint.
Indeed, assume $p(h_ix)=p(h_jy)$ for some $i,j$ and $x,y\in U$, so
that $\gamma h_ix=h_jy$ for some $\gamma\in\Gamma$. Since
$U\subset Y_0$, our assumption on $Y_0$ implies $h^{-1}_j\gamma
h_i\in\Gamma$. But then, since $p$ is injective on $U$, we get
$x=y$, and since the action of $\Gamma$ is free, we conclude that
$h^{-1}_j\gamma h_i=e$. It follows that $i=j$ and $h_ix=h_jy$,
which proves the claim.

Furthermore, the union of the disjoint sets
$p(h_1U),\dots,p(h_nU)$ is the set $\Gamma gZ=p(\Gamma
gp^{-1}(Z))$. Hence, since $N(h_i)=N(g)$ for $i=1,\dots,n$,
$$
\int_{\Gamma g Z}fd\nu=\sum^n_{i=1}\int_{h_iU}f\circ p\,d\mu
=N(g)^{-\beta}\sum^n_{i=1}\int_Uf(p(h_i\cdot))d\mu =N(g)^{-\beta}
R_\Gamma(g)\int_ZT_gfd\nu.
$$

The last assertion of the lemma follows by taking $f=\31_{\Gamma
gZ}$ and observing that then $(T_gf)(z)=1$ for $z\in Z$. \ep

To formulate the next lemma we introduce the following notation.

\begin{definition}
If $\beta\in\R$ and $S$ is a subsemigroup of $G$ containing
$\Gamma$, then we define
$$
\zeta_{S,\Gamma}(\beta):=\sum_{s\in \Gamma\backslash
S}N(s)^{-\beta}=\sum_{s\in\Gamma\backslash
S/\Gamma}N(s)^{-\beta}R_\Gamma(s).
$$
\end{definition}

\begin{lemma} \label{project}
Suppose $\mu$ is as in Proposition~\ref{Scaling1} and $\nu$ is the
measure on $\Gamma\backslash Y$ determined by (\ref{equotientm}).
Assume that the action of $G$ on $X$ is free and that $(G,\Gamma)$
is a Hecke pair. Assume further that $Y_0$ is a measurable
$\Gamma$-invariant subset of $Y$, and $S$ a subsemigroup of $G$
containing $\Gamma$ such that \enu{i} if $gY_0\cap
Y_0\ne\emptyset$ for some $g\in G$ then $g\in\Gamma$; \enu{ii}
$\cup_{s\in S}sY_0$ is a subset of $Y$ of full measure; \enu{iii}
$\zeta_{S,\Gamma}(\beta)<\infty$.

Let $H_S$ be the subspace of $S$-invariant functions in
$L^2(\Gamma\backslash Y,\nu)$, that is, functions $f$ such that
$f(y)=f(sy)$ for all $s\in S$ and a.a. $y\in Y$. Then \enu{1} if
$f\in H_S$ then $\displaystyle
\|f\|^2_2=\zeta_{S,\Gamma}(\beta)\int_{\Gamma\backslash
Y_0}|f(t)|^2d\nu(t)$; \enu{2} the orthogonal projection $P\colon
L^2(\Gamma\backslash Y,d\nu)\to H_S$ is given by
\begin{equation} \label{epro0}
Pf|_{Sy}=\zeta_{S,\Gamma}(\beta)^{-1}\sum_{s\in\Gamma\backslash
S/\Gamma}N(s)^{-\beta}R_\Gamma(s)(T_sf)(y)\ \ \text{for}\ \ y\in
Y_0.
\end{equation}
\end{lemma}

\bp By condition (i) the sets $\Gamma sY_0$ are disjoint for $s$
in different double cosets of $\Gamma$. Since the union of such
sets is the whole space $Y$ (modulo a set of measure zero), by
Lemma~\ref{measure} applied to $Z=\Gamma\backslash Y_0$ for any
$f\in L^2(\Gamma\backslash Y,d\nu)$ we get
\begin{equation} \label{enorm}
\|f\|^2_2=\sum_{s\in\Gamma\backslash S/\Gamma}\int_{\Gamma
sZ}|f|^2d\nu=\sum_{s\in\Gamma\backslash
S/\Gamma}N(s)^{-\beta}R_\Gamma(s)\int_{\Gamma\backslash
Y_0}T_s(|f|^2)d\nu.
\end{equation}
Since $T_s(|f|^2)=|f|^2$ for $f\in H_S$, this gives (1).

\smallskip

Turning to (2), denote by $T$ the operator on
$L^2(\Gamma\backslash Y,d\nu)$ defined by the asserted formula for
$P$. To see that it is well-defined, notice first that the
summation in the right hand side of (\ref{epro0}) is finite for
$f$ in the subspace of $L^2$-functions supported on a finite
collection of sets of the form $p(sY_0)$, $s\in S$, which is a
dense subspace of $L^2(\Gamma\backslash Y,d\nu)$. Thus the
function $Tf$ is well-defined for $f$ in this subspace and,
putting
$\alpha_s=\zeta_{S,\Gamma}(\beta)^{-1}N(s)^{-\beta}R_\Gamma(s)$
and using (\ref{enorm}) twice, we get
$$
\|Tf\|^2_2=\zeta_{S,\Gamma}(\beta)\int_{\Gamma\backslash
Y_0}|Tf|^2d\nu\le \zeta_{S,\Gamma}(\beta)\int_{\Gamma\backslash
Y_0} \left(\sum_{s\in\Gamma\backslash
S/\Gamma}\alpha_sT_s(|f|^2)\right)d\nu=\|f\|^2_2.
$$
It follows that $T$ extends to a well-defined contraction. Since
$Tf=f$ for $f\in H_S$, and the image of~$T$ is~$H_S$, we conclude
that $T=P$. \ep

\bigskip

\section{The Connes-Marcolli system} \label{scm}

Consider the group $G=\glq$ of invertible $2$ by $2$ matrices with
rational coefficients and positive determinant, and its subgroup
$\Gamma=\slz$. For a prime number $p$ consider the field $\Q_p$ of
$p$-adic numbers and its compact subring $\Z_p$ of $p$-adic
integers. We denote by $\af$ the space of finite adeles of $\Q$,
that is, the restricted product of the fields of $\Q_p$ with respect
to $\Z_p$, and by $\Zhat=\prod_p\Z_p$ its maximal compact subring. The
field $\Q$ is a subfield of $\Q_p$, so $\glq$ can be considered as a
subgroup of $\G(\Q_p)$. In particular, we have an action of $\glq$
on $\mtwo(\Q_p)$ by multiplication on the left. Moreover, by
considering the diagonal embedding of $\Q$ into $\af$ we get an
embedding of $\glq$ into $\G(\af)$, and thus an action of $\glq$ on
$\ma$. In addition $\glq$ acts by M\"{o}bius transformations on the
upper halfplane $\HH$. Therefore we have an action of $\glq$ on
$\HH\times\ma$ such that for $g=\fmatr{a}{b}{c}{d}$, $\tau\in\HH$
and $m=(m_p)_p\in\ma$,
$$
g(\tau,(m_p)_p)=\left(\frac{a\tau+b}{c\tau+d},(gm_p)_p\right).
$$
Note that the action of $\slz$ is proper, since already the action
of $\slz$ on $\HH$ is proper.

\smallskip

The $\gl$-system of Connes and Marcolli is now defined as follows,
see~\cite[Section~1.8]{con-mar}.

\begin{definition}
The Connes-Marcolli algebra is the C$^*$-algebra $A=\red{\fibb G
\Gamma Y}$, where $G=\glq$, $\Gamma=\slz$, $G$ acts diagonally on
$X=\HH\times\ma$, and $Y=\HH\times\mtwo(\Zhat)$. The dynamics
$\sigma$ on $A$ is defined by the homomorphism
$N\colon\glq\to\R^*_+$, $N(g)=\det(g)$.
\end{definition}

Notice that since $\Gamma\backslash\HH$ is not compact, the algebra
$A$ is nonunital.

\smallskip

By \cite[Lemma 1.28]{con-mar} the action of $\glq$ on
$X\setminus(\HH\times\{0\})$ is free. Recall briefly the reason.
If $gm=m$ for some prime number $p$ and nonzero $m\in\mtwo(\Q_p)$
then the spectrum of the matrix~$g$ contains $1$, and hence $g$ is
conjugate in $\glq$ to an upper-triangular matrix. But then $g$
has no fixed points in $\HH$. Note that what we have actually
shown is that the action of $\glq$ on
$\HH\times\mtwo(\Q_p)^\times$, where
$\mtwo(\Q_p)^\times=\mtwo(\Q_p)\setminus\{0\}$, is free for any
prime number $p$.

Although the action of $\glq$ on $\HH\times\{0\}$ is not free,
this set can be ignored in the analysis of KMS$_\beta$-states for
$\beta\ne0$, see the proof of \cite[Proposition~1.30]{con-mar}.
Again, recall briefly what happens. Consider the action of $G$ on
$\tilde X=X\setminus (\HH\times\{0\})$, put $\tilde Y=Y\setminus
(\HH\times\{0\})\subset\tilde X$, and then define $I=\red{\fibb G
\Gamma\tilde Y}$. Then $I$ can be considered as an ideal in $A$,
and the quotient algebra $A/I$ is isomorphic to $\red{\fib G
\Gamma \HH}$. Now if $\varphi$ is a $\sigma$-KMS$_\beta$ state on
$A$, the restriction $\varphi|_I$ canonically extends to a
KMS-functional on the multiplier algebra of $I$. Thus we get a
KMS-functional $\tilde\varphi\le\varphi$ on $A$. If
$\tilde\varphi\ne\varphi$ then $\varphi-\tilde\varphi$ is a
positive nonzero KMS-functional on $A$ which vanishes on $I$.
Hence we get a KMS-state on $A/I\cong\red{\fib G \Gamma \HH}$. By
Lemma~\ref{HeckeM} the multiplier algebra of $\red{\fib G \Gamma
\HH}$ contains the reduced Hecke C$^*$-algebra
$\redheck{G}{\Gamma}$. The latter algebra contains in turn the
C$^*$-algebra of $Z(G)/(Z(G)\cap\Gamma)$, where $Z(G)$ is the
center of $\glq$, that is, the group of scalar matrices. But since
the dynamics scales nontrivially some unitaries in this algebra,
the algebra can not have any KMS$_\beta$-states for $\beta\ne0$.
This contradiction shows that $\varphi=\tilde\varphi$, so that
$\varphi$ is completely determined by~$\varphi|_I$.

Since the action of $G$ on $\tilde X=\HH\times\ma^\times$, where
$\ma^\times=\ma\setminus\{0\}$, is free, we can apply
Proposition~\ref{Scaling1} and conclude that there is a one-to-one
correspondence between KMS$_\beta$-weights on $I$ with domain of
definition containing $C_c(\Gamma\backslash\tilde Y)$ and measures
$\mu$ on $\tilde Y=\HH\times\mr^\times$ such that
$\mu(gZ)=\det(g)^{-\beta}\mu(Z)$ if both $Z$ and $gZ$ are subsets of
$\tilde Y$. By Lemma~\ref{extensionX} we can uniquely extend any
such measure to a measure on $\tilde X=G\tilde
Y=\HH\times\ma^\times$ such that $\mu(gZ)=\det(g)^{-\beta}\mu(Z)$ for
$Z\subset\tilde X$. To get a state on
$I=\red{\fibb{G}{\Gamma}{\tilde Y}}$ we need the normalization
condition $\mu(\Gamma\backslash \tilde Y)=1$ (that is, the
$\Gamma$-invariant measure $\mu$ on $\tilde Y$ defines a probability
measure on $\Gamma\backslash\tilde Y$). Note also that if
$\beta\ne0$ and we have a measure on $X=\HH\times\ma$ with the same
properties as above, then $\HH\times\ma^\times$ is a subset of full
measure, since scalar matrices act trivially on $\HH$  and so $\HH$
cannot support a measure scaled nontrivially by them.

\smallskip

Summarizing the above discussion we get the following.

\begin{proposition} \label{cm1}
For $\beta\ne0$ there is a one-to-one correspondence between
$\sigma$-KMS$_\beta$-states on the Connes-Marcolli system and
$\Gamma$-invariant measures $\mu$ on $\HH\times\ma$ such that
$$
\mu(\Gamma\backslash(\HH\times\mtwo(\Zhat)))=1\ \ \text{and}\ \
\mu(gZ)=\det(g)^{-\beta}\mu(Z)
$$
for any $g\in\glq$ and compact $Z\subset \HH\times\ma$.
\end{proposition}

Denote by $\mtwo^i(\af)$ the set of matrices $m=(m_p)_p\in\ma$
such that $\det(m_p)\ne0$ for every prime $p$. Notice that
$\mtwo^i(\af)$ is the set of non zero-divisors in $\ma$. Our next
goal is to show that if $\beta\ne0,1$ then $\HH\times\mtwo^i(\af)$
is a subset of full measure for any measure $\mu$ as in
Proposition~\ref{cm1}.

First let us recall the following simple properties of the Hecke
pair $(G,\Gamma)=(\glq,\slz)$. Put $\mzp=\glq\cap\mz$.

\begin{lemma} \label{class}
Every double coset of $\Gamma$ in $\mzp$ has a unique
representative of the form $\diag a d$ with $a,d\in\N$ and $a|d$.
Furthermore,
$$
R_\Gamma\diag{a}{d}=\frac{d}{a}\prod_{p\ \text{prime}\colon
pa|d}(1+p^{-1}),
$$
and as representatives of the right cosets of $\Gamma$ contained
in $\Gamma{\diag a d}\Gamma$ we can take the matrices
$$
\fmatr{ak}{am}{0}{al}
$$
with $k,l\in \N$ and $m\in\Z$ such that $kl=d/a$, $0\le m<l$ and
$\gcd(k,l,m)=1$.

In particular, $R_\Gamma(g)=R_\Gamma(g^{-1})$ for every
$g\in\glq$.
\end{lemma}

\bp See e.g.~\cite[Chapter~IV]{kri}. \ep

For a prime $p$ put $G_p=\glp(\Z[p^{-1}])\subset\glq$. Observe that
if $g\in G_p$ then $\det(g)$ is a power of~$p$, and if we multiply
$g$ by a sufficiently large power of $\diag p p$, we get an element
in $\mzp$ with determinant a power of $p$. But by Lemma~\ref{class}
the double coset of $\Gamma$ containing such an element has a
representative of the form $\diag{p^k}{p^l}$, $0\le k\le l$. We may
therefore conclude that $G_p$ is the subgroup of $\glq$ generated by
$\Gamma$ and $\diag 1 p$. Using that a positive rational number is a
power of $p$ if and only if it belongs to the group of units
$\Z_q^*$ of the ring $\Z_q$ for all primes $q\ne p$, we may also
conclude that $g\in\glq$ belongs to $G_p$ if and only if it belongs
to $\G(\Z_q)$ for all $q\ne p$.

\begin{lemma} \label{decomp}
We have $\G(\Q_p)=G_p\G(\Z_p)$.
\end{lemma}

\bp Let $r\in\G(\Q_p)$. Then $r\Z^2_p$ is a $\Z_p$-lattice in
$\Q_p^2$, that is, an open compact $\Z_p$-submodule.
By~\cite[Theorem~V.2]{weil} there exists a subgroup $L\cong\Z^2$
of $\Q^2$ such that the closure of $L$ in $\Q_p^2$ coincides
with~$r\Z^2_p$, and the closure of $L$ in $\Q^2_q$ is $\Z^2_q$ for
$q\ne p$. Choose $g\in\glq$ such that $g\Z^2=L$. Since
$g\Z^2_p=r\Z^2_p$, we have $g^{-1}r\in\G(\Z_p)$. Since
$g\Z^2_q=\Z^2_q$ for $q\ne p$, we also have $g\in\G(\Z_q)$. Hence
$g\in G_p$. \ep

It is also possible to give an elementary proof of
Lemma~\ref{decomp} using matrix factorization and density of
$\Z[p^{-1}]$ in $\Q_p$.

\begin{lemma} \label{singular}
Let $p$ be a prime and $\mu_p$ a $\Gamma$-invariant measure on
$\HH\times\mtwo(\Q_p)$ such that
$$
\mu_p(\HH\times\{0\})=0,\ \
\mu_p(\Gamma\backslash(\HH\times\mtwo(\Z_p)))<\infty\ \ \text{and}\ \
\mu_p(gZ)=\det(g)^{-\beta}\mu_p(Z)
$$
for $g\in G_p$ and $Z\subset \HH\times\mtwo(\Q_p)$. If
$\beta\ne1$, then the set $\HH\times\G(\Q_p)$ is a subset of full
measure in~$\HH\times\mtwo(\Q_p)$.
\end{lemma}

\bp Denote by $\tilde\nu$ the measure on
$\Gamma\backslash(\HH\times\mtwo(\Q_p))$ defined by the
$\Gamma$-invariant measure $\mu_p$. For a $\Gamma$-invariant subset
$Z$ of $\mtwo(\Q_p)$, the set $\HH\times Z$ is $\Gamma$-invariant.
We can thus define a measure $\nu$ on the $\sigma$-algebra of
$\Gamma$-invariant Borel subsets of $\mtwo(\Q_p)$ by
$\nu(Z)=\tilde\nu(\Gamma\backslash(\HH\times Z))$. Note that since
the action of $\Gamma$ on $\mtwo(\Q_p)$ is not proper and,
accordingly, the quotient space $\Gamma\backslash \mtwo(\Q_p)$ is
quite bad, we do not want to consider $\Gamma$-invariant subsets
of $\mtwo(\Q_p)$ as subsets of this quotient space, and do not try
to define a measure on all Borel subsets of $\mtwo(\Q_p)$ out of
$\nu$.

If $g\in G_p$ and $f$ is a positive Borel $\Gamma$-invariant
function on $\mtwo(\Q_p)$ then by Lemma~\ref{scaling2} applied to
the function $F\colon(\tau,m)\mapsto f(m)$ on
$\Gamma\backslash(\HH\times\mtwo(\Q_p))$ we conclude that
\begin{equation}\label{ehecksing}
\begin{split}
\int_{\mtwo(\Q_p)} T_gfd\nu&=
\int_{\Gamma\backslash(\HH\times\mtwo(\Q_p))} T_gFd\tilde\nu\\ &=
\det(g)^\beta\int_{\Gamma\backslash(\HH\times\mtwo(\Q_p))}
Fd\tilde\nu=\det(g)^\beta\int_{\mtwo(\Q_p)} fd\nu.
\end{split}
\end{equation}
By assumption we also have $\nu(\mtwo(\Z_p))<\infty$. We have to
show that the measure of the set of nonzero singular matrices is
zero.

We claim that the set of nonzero singular matrices with
coefficients in $\Q_p$ is the disjoint union of the sets
$$
Z_k=\sltwo(\Z_p)\diag{0}{p^k}\G(\Z_p),\ k\in\Z.
$$
This is proved in a standard way: given a nonzero singular matrix
we use multiplication by elements of $\G(\Z_p)$ on the right to
get a matrix with zero first column, and then multiplication by
elements of $\sltwo(\Z_p)$ on the left to get the required form.
To show that the sets do not intersect, observe that the maximum
of the $p$-adic valuations of the coefficients of a matrix does
not change under multiplication by elements of $\G(\Z_p)$ on
either side.

Consider the functions $f_k=\31_{Z_k}$, $k\in\Z$. For
$g=\diag{1}{p^{-1}}$ we claim that
$$
T_gf_0=\frac{1}{p+1}f_0+\frac{p}{p+1}f_1.
$$
Indeed, since the action of $G_p$ commutes with the right action of
$\G(\Z_p)$, the function $T_gf_0$ is $\G(\Z_p)$-invariant. On the
other hand, the sets $Z_k$ are clopen subsets of the set of singular
matrices, so that the function $f_0$ is continuous on this set. But
then $T_gf_0$ is also continuous. Since $T_gf_0$ is
$\Gamma$-invariant, and $\Gamma$ is dense in $\sltwo(\Z_p)$ (see
e.g. \cite[Lemma~1.38]{shi} for an elementary proof of a stronger
result: $\Gamma$ is dense in $\slr$), we conclude that $T_gf_0$ is
left $\sltwo(\Z_p)$-invariant. Hence $T_gf_0$ is constant on the
sets $Z_k$. So to prove the above identity it suffices to check it
on the matrices $\diag{0}{p^k}$. Since
$g=\diag{p^{-1}}{p^{-1}}\diag{p}{1}$, by Lemma~\ref{class} we can
take the matrices
$$
\diag{1}{p^{-1}},\ \ \fmatr{p^{-1}}{np^{-1}}{0}{1}, \ 0\le n\le
p-1,
$$
as representatives of the right cosets of $\Gamma$ contained in
$\Gamma g\Gamma$. Then
$$
(T_gf_0)\diag{0}{p^k}=\frac{1}{p+1}f_0
\diag{0}{p^{k-1}}+\frac{1}{p+1}\sum^{p-1}_{n=0}
f_0\fmatr{0}{np^{k-1}}{0}{p^k}.
$$
Since the matrices $\diag{0}{p^{k-1}}$ and
$\fmatr{0}{np^{k-1}}{0}{p^k}$, $1\le n\le p-1$, belong to
$Z_{k-1}$, we see that
$$
T_gf_0|_{Z_1}=\frac{p}{p+1},\ \ T_gf_0|_{Z_0}=\frac{1}{p+1},\ \
T_gf_0|_{Z_k}=0\ \ \text{for}\ \ k\ne0,1,
$$
and this is exactly what we claimed.

It follows from \eqref{ehecksing} that
$$
p^{-\beta}\nu(Z_0)=\frac{1}{p+1}\nu(Z_0)+\frac{p}{p+1}\nu(Z_1).
$$
On the other hand, for $g=\diag{p^{-1}}{p^{-1}}$ we get
$T_gf_k=f_{k+1}$, so that
$$
p^{-2\beta}\nu(Z_k)=\nu(Z_{k+1}).
$$
If $\nu(Z_0)\ne0$ this implies that $p^{-\beta}$ is a solution of
the quadratic equation
$$
(p+1)x=1+p\,x^2,
$$
Thus either $p^{-\beta}=p^{-1}$ or $p^{-\beta}=1$. Since
$\beta\ne1$ we get $\beta=0$. But then $\nu(Z_k)=\nu(Z_0)$ for any
$k$, and this contradicts $\nu(\mtwo(\Z_p))<\infty$. The
contradiction shows that $\nu(Z_0)=0$. Hence $\nu(Z_k)=0$ for any
$k$, and we conclude that the measure of the set of singular
matrices is zero. \ep

We are now ready to show that for $\beta\ne0,1$ the set
$\ma\setminus\mtwo^i(\af)$ of zero-divisors has measure zero.

\begin{corollary} \label{singular2}
Assume $\beta\ne0,1$ and $\mu$ is a measure with properties as in
Proposition~\ref{cm1}. Then $\HH\times\mtwo^i(\af)$ is a subset of
full measure in $\HH\times\ma$.
\end{corollary}

\bp Fix a prime $p$. First of all note that the set
$$
\{(\tau,m)\in\HH\times\ma\mid m_p=0\}
$$
has measure zero. Indeed, as we already remarked before
Proposition~\ref{cm1}, the set $\HH\times\{0\}$ has measure zero.
So if our claim is not true, the set
$$
\{(\tau,m)\in\HH\times\mtwo(\Zhat)^\times\mid m_p=0\}
$$
has positive measure. Since the action of $\Gamma$ on this set is
free, there is a subset $U$ of positive measure such that $\gamma
U\cap U=\emptyset$ for $\gamma\in\Gamma$, $\gamma\ne e$. Then for
$g=\diag{p}{p}$ the set $U_k=g^kU$, $k\in\Z$, still has the property
that $\gamma U_k\cap U_k=\emptyset$ for $\gamma\in\Gamma$,
$\gamma\ne e$, since $g$ commutes with $\Gamma$. As $U_k$ is
contained in $\HH\times\mtwo(\Zhat)$, it follows that $\mu(U_k)\le1$.
On the other hand, $\mu(U_k)=p^{-2\beta k}\mu(U)$. Letting
$k\to-\infty$ if $\beta>0$ and $k\to+\infty$ if $\beta<0$, we get a
contradiction.

Consider now the restriction of $\mu$ to the set
$$
\HH\times\mtwo(\Q_p)\times\prod_{q\ne p}\mtwo(\Z_q),
$$
and use the projection onto the first two factors to get a measure
$\mu_p$ on $\HH\times\mtwo(\Q_p)$. By the first part of the proof
the set $\HH\times\{0\}$ has $\mu_p$-measure zero. Since the image
of $G_p$ in $\G(\Q_q)$ lies in $\G(\Z_q)$ for $q\ne p$, the
scaling property of $\mu$ implies that
$$
\mu_p(gZ)=\det(g)^{-\beta}\mu_p(Z)\ \ \text{for}\ \ Z\subset
\HH\times\mtwo(\Q_p),\ g\in G_p.
$$
Since the action of $\Gamma$ on $\HH\times\mtwo(\Q_p)^\times$ is
free, the normalization condition on $\mu$ implies that
$\mu_p(\Gamma\backslash(\HH\times\mtwo(\Z_p))=1$. Thus $\mu_p$
satisfies the assumptions of Lemma~\ref{singular}. Hence
$\HH\times\G(\Q_p)$ is a set of full $\mu_p$-measure. This means
that the set of points $(\tau,m)\in\HH\times\mtwo(\Zhat)$ with
$\det(m_p)=0$ has $\mu$-measure zero. By taking the union of such
sets for all primes $p$ and multiplying it by elements of~$\glq$ we
get a set of measure zero, which is the complement of the set
$\HH\times\mtwo^i(\af)$. \ep

To get further properties of a measure $\mu$ as above, let us
recall the following well-known computation, see e.g.
\cite[Section 3.2]{shi} for more general results on formal
Dirichlet series. Denote by $S_p$ the semigroup $G_p\cap\mzp$.
Alternatively, $S_p$ is the set of elements $m\in\mzp$ with
determinant a nonnegative power of $p$. Then from
Lemma~\ref{class} we know that as representatives of the right
cosets of $\Gamma$ in $S_p$ we can take the matrices
$\fmatr{p^k}{m}{0}{p^l}$, $k,l\ge0$, $0\le m<p^l$. Therefore
\begin{equation} \label{ezeta}
\zeta_{S_p,\Gamma}(\beta)=\sum_{s\in\Gamma\backslash
S_p}\det(s)^{-\beta}=\sum^{\infty}_{k,l=0}p^{-\beta(k+l)}p^l
=\begin{cases}+\infty,& \text{if}\ \beta\le1,\cr
(1-p^{-\beta})^{-1}(1-p^{-\beta+1})^{-1}, & \text{if}\
\beta>1.\end{cases}
\end{equation}

Since $\Gamma=G_p\cap\G(\Z_p)$, we can apply Lemma~\ref{measure} to
the group $G_p$ acting on $\HH\times\ma^\times$ and the set
$$
Y_0=\HH\times\G(\Z_p)\times\prod_{q\ne p}\mtwo(\Z_q).
$$
Then for any $s\in S_p$ we get
$$
\mu(\Gamma\backslash\Gamma sY_0)=\det(s)^{-\beta}R_\Gamma(s)
\mu(\Gamma\backslash Y_0).
$$
The sets $\Gamma s Y_0$ are disjoint for $s$ in different double
cosets of $\Gamma$, and their union is the set
$$
\HH\times\mtwo^i(\Z_p)\times\prod_{q\ne p}\mtwo(\Z_q),
$$
where $\mtwo^i(\Z_p)=\mtwo(\Z_p)\cap\G(\Q_p)$. By
Corollary~\ref{singular2} the above set is a subset of
$\HH\times\mtwo(\Zhat)$ of full measure for $\beta\ne0,1$. Therefore
we obtain
\begin{equation} \label{emass}
1=\sum_{s\in\Gamma\backslash S_p/\Gamma}\mu(\Gamma\backslash\Gamma
s Y_0)=\sum_{s\in\Gamma\backslash
S_p/\Gamma}\det(s)^{-\beta}R_\Gamma(s)\mu(\Gamma\backslash Y_0)
=\zeta_{S_p,\Gamma}(\beta)\mu(\Gamma\backslash Y_0).
\end{equation}
This gives a contradiction if $\beta<1$. Thus for $\beta<1$,
$\beta\ne0$, there are no KMS$_\beta$-states. On the other hand,
for $\beta>1$ we get
$$
\mu(\Gamma\backslash Y_0)=\zeta_{S_p,\Gamma}(\beta)^{-1}
=(1-p^{-\beta})(1-p^{-\beta+1}).
$$

Assuming now that $\beta>1$ we can perform a similar computation for
any finite set of primes instead of just one prime. Given a finite
set $F$ of primes consider the group $G_F$ generated by $G_p$ for
all $p\in F$. Put also $S_F=\mzp\cap G_F$. Then $S_F$ is the set of
matrices $m\in\mzp$ such that all prime divisors of $\det(m)$ belong
to $F$. Let
$$
Y_F=\HH\times\left(\prod_{p\in F}\G(\Z_p)\right)\times
\left(\prod_{q\notin F}\mtwo(\Z_q)\right).
$$
Then a computation similar to (\ref{ezeta}) and (\ref{emass})
yields
\begin{equation} \label{ezeta1}
\zeta_{S_F,\Gamma}(\beta)=\prod_{p\in F}(1-p^{-\beta})^{-1}
(1-p^{-\beta+1})^{-1}\ \ \text{and}\ \ \mu(\Gamma\backslash
Y_F)=\prod_{p\in F}(1-p^{-\beta})(1-p^{-\beta+1}).
\end{equation}
The intersection of the sets $Y_F$ over all finite subsets $F$ of
prime numbers is the set $\HH\times\G(\Zhat)$. So for $\beta>2$ we
get
$$
\mu(\Gamma\backslash(\HH\times\G(\Zhat)))=\prod_p
(1-p^{-\beta})(1-p^{-\beta+1})=\zeta(\beta)^{-1}\zeta(\beta-1)^{-1},
$$
where $\zeta$ is the Riemann $\zeta$-function. On the other hand,
for $\beta\in(1,2]$ we get
$\mu(\Gamma\backslash(\HH\times\G(\Zhat)))=0$.

Assume now that $\beta>2$. In this case similarly to (\ref{ezeta})
we have
$$
\zeta_{\mzp,\Gamma}(\beta)=\zeta(\beta)\zeta(\beta-1).
$$
So analogously to (\ref{emass}) we get
$$
\mu(\Gamma\backslash\mzp(\HH\times\G(\Zhat)))
=\zeta_{\mzp,\Gamma}(\beta)\mu(\Gamma\backslash(\HH\times\G(\Zhat)))
=1.
$$
We thus see that $\mzp(\HH\times\G(\Zhat))$ is  a subset of
$\HH\times\mtwo(\Zhat)$ of full measure. Hence
$\glq(\HH\times\G(\Zhat))$ is a subset of $\HH\times\ma$ of full
measure. By Lemma~\ref{decomp} the set $\glq(\HH\times\G(\Zhat))$ is
nothing but $\HH\times\G(\af)$.

To summarize, we have shown that for $\beta>2$ the problem of
finding all measures $\mu$ on $\HH\times\ma$ satisfying the
conditions in Proposition~\ref{cm1} reduces to finding all
measures on $\HH\times\G(\af)$ such that
$$
\mu(gZ)=\det(g)^{-\beta}\mu(Z)\ \ \text{and}\ \
\mu(\Gamma\backslash(\HH\times\G(\Zhat)))=\zeta(\beta)^{-1}
\zeta(\beta-1)^{-1}.
$$
By Lemma~\ref{extension} any $\Gamma$-invariant measure on
$\HH\times\G(\Zhat)$ extends uniquely to a measure on
$\HH\times\G(\af)$ satisfying the scaling condition. Thus we get a
one-to-one correspondence between measures $\mu$ with properties
as in Proposition~\ref{cm1} and measures on
$\Gamma\backslash(\HH\times\G(\Zhat))$ of total mass
$\zeta(\beta)^{-1}\zeta(\beta-1)^{-1}$. Clearly, extremal measures
$\mu$ correspond to point masses.

\smallskip

We have thus recovered the following result of Connes and
Marcolli~\cite[Theorem~1.26 and Corollary~1.32]{con-mar}.

\begin{theorem} \label{conmar}
For the Connes-Marcolli $\G$-system we have: \enu{i} for
$\beta\in(-\infty,0)\cup(0,1)$ there are no KMS$_\beta$-states;
\enu{ii} for $\beta>2$ there is a one-to-one affine correspondence
between KMS$_\beta$-states and probability measures on
$\Gamma\backslash(\HH\times\G(\Zhat))$; in particular, extremal
KMS$_\beta$-states are in bijection with $\Gamma$-orbits in
$\HH\times\G(\Zhat)$.
\end{theorem}

\begin{remark} This is not exactly what is stated
in~\cite{con-mar}. First of all, the cases $\beta=0,1$ require
considerations with singular matrices, and in these cases we do have
KMS-states, see Remark~\ref{rsingular} below. Secondly, the
classification of extremal KMS$_\beta$-states for $\beta>2$
in~\cite[Theorem~1.26]{con-mar} is in terms of invertible
$\Q$-lattices up to scaling. To see that our
Theorem~\ref{conmar}(ii) says the same, recall that the isomorphism
from~\cite[Equation~(1.87)]{con-mar} identifies
$\Gamma\backslash(\HH\times\G(\Zhat))$ with the set of invertible
$\Q$-lattices in~$\C$ up to scaling, and observe that the state
$\varphi_{\beta,l}$ defined in~\cite[Theorem~1.26(ii)]{con-mar}
associated with $l=(\tau,\rho)\in\HH\times\G(\Zhat)$ is exactly the
KMS$_\beta$-state corresponding to the orbit $\Gamma(\tau,\rho)$.
Since the $\Q$-lattice picture will not be used in the remaining
part of the paper, we omit the details.
\end{remark}

\bigskip

\section{Uniqueness of the KMS$_\beta$-state in the critical region
$1<\beta\le2$} \label{scm2}

In this section we analyze the Connes-Marcolli system in the
region $\beta\in(1,2]$.

\smallskip

For each such $\beta$ let us first construct a KMS$_\beta$-state,
or equivalently, a measure $\mu_\beta$ on $\HH\times\ma$
satisfying the conditions in Proposition~\ref{cm1}.

For each prime number $p$ consider the Haar measure on $\G(\Z_p)$
normalized such that the total mass is
$(1-p^{-\beta})(1-p^{-\beta+1})$. By the same argument as in the
proof of Lemma~\ref{extension}, this measure extends to a unique
measure $\mu_{\beta,p}$ on $\G(\Q_p)$ such that
$$
\mu_{\beta,p}(Z)=\sum_{g\in\G(\Z_p)\backslash\G(\Q_p)}
|\det(g)|_p^{-\beta}\mu_{\beta,p}(gZ\cap\G(\Z_p))
$$
for compact $Z\subset\G(\Q_p)$, where $|a|_p$ denotes the $p$-adic
valuation of $a$. The measure $\mu_{\beta,p}$ satisfies
$$
\mu_{\beta,p}(gZ)=|\det(g)|^\beta_p\mu_{\beta,p}(Z)\ \ \text{for}\ \
g\in\G(\Q_p).
$$
Since $|\det(g)|_p=1$ for $g\in\G(\Z_p)$, it is clear that
$\mu_{\beta,p}$ is left $\G(\Z_p)$-invariant. But since the Haar
measure on $\G(\Z_p)$ is biinvariant, we conclude that
$\mu_{\beta,p}$ is also right $\G(\Z_p)$-invariant. By setting
$\mu_{\beta,p}(Z)=\mu_{\beta, p}(Z\cap\G(\Q_p))$ for Borel $Z\subset
\mtwo(\Q_p)$ we extend $\mu_{\beta, p}$ to a measure
on~$\mtwo(\Q_p)$. Using that $\mtwo^i(\Z_p)=S_p\G(\Z_p)$, similarly
to (\ref{emass}) we find
$$
\mu_{\beta,p}(\mtwo(\Z_p))=\zeta_{S_p,\Gamma}(\beta)
\mu_{\beta,p}(\G(\Z_p))=1.
$$
Hence we can define a measure on $\ma$ by
$\mu_{\beta,f}=\prod_p\mu_{\beta,p}$. By construction and
Lemma~\ref{singular} this is the unique product-measure such that
$\mu_{\beta,f}(\mtwo(\Zhat))=1$ and
\begin{equation} \label{escalingfull}
\mu_{\beta,f}(gZr)=\left(\prod_p|\det(g_p)|_p\right)^\beta\mu_{\beta,f}(Z)
\end{equation}
for $Z\subset\ma$, $g=(g_p)_p\in\G(\af)$ and $r\in\G(\Zhat)$.  Note
that since a Haar measure on the additive group $\ma$ is a
product-measure satisfying (\ref{escalingfull}) with $\beta=2$, we
see that  $\mu_{2,f}$ is a Haar measure on~$\ma$.

Denote by $\mu_\infty$ the unique $\glq$-invariant measure on $\HH$
such that $\mu_\infty(\Gamma\backslash\HH)=1$.

Now put $\mu_\beta=2\mu_\infty\times\mu_{\beta,f}$. Then $\mu_\beta$
satisfies the conditions in Proposition~\ref{cm1}, so it corresponds
to a KMS$_\beta$-state on the Connes-Marcolli C$^*$-algebra. Indeed,
the scaling condition is satisfied since $\prod_p|q|_p=q^{-1}$ for
$q\in\Q^*_+$. The factor $2$ is needed for the normalization
condition, since the element $-1\in\Gamma$ acts trivially on~$\HH$,
while $\mu_{\beta,f}(\{\pm1\}\backslash\mtwo(\Zhat))=1/2$.

Note that the construction of $\mu_\beta$ makes sense for all
$\beta>1$.

\smallskip

We can now formulate our main result.

\begin{theorem} \label{cm2}
For each $\beta\in(1,2]$ the state corresponding to the measure
$\mu_\beta$ is the unique KMS$_\beta$-state on the Connes-Marcolli
system.
\end{theorem}

We shall prove a slightly stronger result which may look more
natural if one leaves aside the motivation for the Connes-Marcolli
system. Namely, we replace $\HH$ by $\pgl=\glp(\R)/\R^*$. Recall
that $\pgl$ acts transitively on $\HH$, and $\so/\{\pm1\}$ is the
stabilizer of the point $i\in\HH$, so that $\HH=\pgl/\pso$. Denote by
$\bar\mu_\infty$ the Haar measure on $\pgl$ normalized such that
$\bar\mu_\infty(\Gamma\backslash\pgl)=1$. Define then a measure on
$\pgl\times\ma$ by
$\bar\mu_\beta=2\bar\mu_\infty\times\mu_{\beta,f}$.

\begin{theorem} \label{pgl}
For $\beta\in(1,2]$ the measure $\bar\mu_\beta$ is the unique
$\Gamma$-invariant measure on the space $\pgl\times\ma$ such that
$$
\bar\mu_\beta(\Gamma\backslash(\pgl\times\mtwo(\Zhat)))=1\ \
\text{and}\ \ \bar\mu_\beta(gZ)=\det(g)^{-\beta}\bar\mu_\beta(Z)
$$
for compact $Z\subset\pgl\times\ma$ and $g\in\glq$.
\end{theorem}

Theorem~\ref{cm2} follows from the above theorem since every measure
$\mu$ on $\HH\times\ma$ satisfying the conditions in
Proposition~\ref{cm1} gives rise to a measure $\bar\mu$ on
$\pgl\times\ma$ satisfying the conditions in Theorem~\ref{pgl} by
the formula
$$
\int_{\pgl\times\ma}fd\bar\mu =\int_{\HH\times\ma}\left(\int_{\pso}
f(\cdot\,g)dg\right)d\mu,
$$
where for $x=(h,m)\in\pgl\times\ma$ and $g\in\pso$ we put
$xg=(hg,m)$, and different measures $\mu$ give rise to different
$\bar\mu$'s.

\smallskip

Turning to the proof of Theorem~\ref{pgl} our first goal is to
show uniqueness of $\bar\mu_\beta$ under the additional assumption
of invariance under the right action of $\G(\Zhat)$ on $\ma$.

\smallskip

Let $F$ be a finite set of prime numbers. Recall that we denote by
$S_F$ the semigroup of matrices $m\in\mzp$ such that all prime
divisors of $\det(m)$ belong to $F$. We then introduce an
operator~$T_F$ on the space of bounded functions on
$\Gamma\backslash \pgl$ by
\begin{equation} \label{edefpro}
(T_Ff)(\tau)=\zeta_{S_F,\Gamma}(\beta)^{-1}\sum_{s\in
\Gamma\backslash S_F/\Gamma}\det(s)^{-\beta}R_\Gamma(s)(T_sf)(\tau).
\end{equation}
Denote by $\bar\nu_\infty$ the measure on $\Gamma\backslash\pgl$
defined by $\bar\mu_\infty$. The following result is a key point in
our argument for uniqueness of the $\G(\Zhat)$-invariant measure.

\begin{lemma} \label{clozel}
For any finite set $J$ of prime numbers, $f\in
C_c(\Gamma\backslash\pgl)$, $\eps>0$ and compact subset
$\Omega\subset\Gamma\backslash\pgl$, there exists a finite set $F$
of prime numbers that is disjoint from $J$ and satisfies
$$
\left|(T_Ff)(\tau)-\int_{\Gamma\backslash\pgl}f
d\bar\nu_\infty\right|<\eps\ \ \text{for all}\ \ \tau\in\Omega.
$$
\end{lemma}

\bp By~\cite[Theorem~1.7]{cou} and Remark~(3) following it, see
also~\cite{eso} for an alternative proof of a slightly weaker
result, there exists a constant $M$ such that
$$
\left|(T_gf)(\tau)-\int_{\Gamma\backslash\pgl}f
d\bar\nu_\infty\right|< \frac{\eps}{2}
$$
for $\tau\in\Omega$ and any $g\in\glq$ with $R_\Gamma(g)>M$. We may
assume that $M$ is such that $p<M$ for any $p\in J$. Let $F$ be a
finite set of prime numbers greater than $M$. Then from
Lemma~\ref{class} we see that $R_\Gamma(s)>M$ for any $s\in S_F$
such that $\Gamma s\Gamma$ contains a nonscalar diagonal matrix. On
the other hand,
$$
\sum_{\substack{s\in \Gamma\backslash S_F/\Gamma:\\s\
\text{scalar}}}\det(s)^{-\beta}=\prod_{p\in
F}\left(\sum^\infty_{k=0}p^{-2\beta k}\right)=\prod_{p\in
F}(1-p^{-2\beta})^{-1}\le\zeta(2\beta).
$$
Since the operators $T_g$ are contractions in the supremum-norm,
we can find $C>0$ such that
$$
\left|(T_gf)(\tau)-\int_{\Gamma\backslash\pgl}f
d\bar\nu_\infty\right|\le C\ \ \text{for}\ \ \tau\in\Omega\ \
\text{and}\ \ g\in\glq.
$$
Therefore by considering separately the summation over double
cosets with nonscalar and scalar representatives we get
$$
\left|(T_Ff)(\tau)-\int_{\Gamma\backslash\pgl}f
d\bar\nu_\infty\right|\le
\frac{\eps}{2}+\frac{\zeta(2\beta)}{\zeta_{S_F,\Gamma}(\beta)}C\ \
\text{for any}\ \ \tau\in\Omega.
$$
Recall that by (\ref{ezeta1})
$$
\zeta_{S_F,\Gamma}(\beta)=\prod_{p\in F}(1-p^{-\beta})^{-1}
(1-p^{-\beta+1})^{-1}.
$$
Since for $\beta\le2$ this product diverges as $F$ increases, we see
that by choosing sufficiently large~$F$ we can make the second
summand in the estimate above arbitrarily small, hence we are done.
\ep

We can now analyze the case of measures on $\pgl\times\ma$ that are
invariant under the right action of $\G(\Zhat)$ on the second factor.

\begin{lemma} \label{symm}
The measure $\bar\mu_\beta$ is the unique right
$\G(\Zhat)$-invariant measure on $\pgl\times\ma$ that satisfies the
conditions in Theorem~\ref{pgl}. Furthermore, the action of $\glq$
on the space $(\pgl\times(\ma/\G(\Zhat)),\bar\mu_\beta)$ is ergodic.
\end{lemma}

\bp The measure $\bar\mu_\beta$ is right $\G(\Zhat)$-invariant and
satisfies the conditions in Theorem~\ref{pgl} by construction.
Suppose $\tilde\mu$ is another such measure. Let $\tilde\nu$ and
$\bar\nu_\beta$ be the measures on the quotient space
$\Gamma\backslash(\pgl\times\ma)$ defined by $\tilde\mu$ and
$\bar\mu_\beta$, respectively. Let $H$ be the subspace of
$\mzp$-invariant functions in
$L^2(\Gamma\backslash(\pgl\times\mtwo(\Zhat)),d\tilde\nu)$, and denote
by $P$ the orthogonal projection onto $H$. Our first goal is to
compute how $P$ acts on $\G(\Zhat)$-invariant functions.

Let $F$ be a nonempty finite set of prime numbers. Apply
Lemma~\ref{project}(2) to the group $G_F$, the semigroup~$S_F$, the
set $Y=\pgl\times\mtwo(\Zhat)$ and the subset
$$
Y_F=\pgl\times\prod_{p\in F}\G(\Z_p) \times\prod_{q\notin
F}\mtwo(\Z_q)
$$
in place of $Y_0$. Note that we can do this because $S_FY_F$
coincides with
$$
\pgl\times\prod_{p\in F}\mtwo^i(\Z_p)\times \prod_{q\notin
F}\mtwo(\Z_q),
$$
which by Corollary~\ref{singular2} (or rather its analogue with
$\HH$ replaced by $\pgl$) is a subset of $Y$ of full measure. Thus,
denoting by $P_F$ the projection onto the subspace of
$S_F$-invariant functions, for $f_0\in
L^2(\Gamma\backslash(\pgl\times\mtwo(\Zhat)),d\tilde\nu)$ we have
\begin{equation} \label{epro}
P_Ff_0|_{S_Fx}=\zeta_{S_F,\Gamma}(\beta)^{-1} \sum_{s\in
\Gamma\backslash S_F/\Gamma}\det(s)^{-\beta}R_\Gamma(s)(T_sf_0)(x)\
\ \text{for each} \ \ x\in Y_F.
\end{equation}
Given a finite set $J$ of prime numbers which is disjoint from $F$,
and a bounded Borel function $f$ on $\Gamma\backslash\pgl$, apply
\eqref{epro} to the function $f_0=f_J$, where $f_J$ is defined by
$$
f_J(x)=\begin{cases}
f(\tau), & \text{if}\ x=(\tau,m)\in Y_J,\\
0, &\text{otherwise.}\end{cases}
$$
Then using the operator $T_F$ defined in (\ref{edefpro}), we can
write
$$
P_Ff_J=(T_Ff)_J.
$$
Assume now that $f$ is continuous and compactly supported.
By Lemma~\ref{clozel} we can find a sequence $\{F_n\}_n$ of finite
sets disjoint from $J$ such that $\{T_{F_n}f\}_n$ converges to $\int
fd\bar\nu_\infty$ uniformly on compact sets. Hence the sequence
$\{P_{F_n}f_J\}_n$ converges weakly in~$L^2$ to $\int
fd\bar\nu_\infty\,(\31_{\Gamma\backslash \pgl})_J= \int
fd\bar\nu_\infty\,\31_{\Gamma\backslash Y_J}$. Since $PP_F=P$ for
every $F$, we get
$$
Pf_J=\int fd\bar\nu_\infty\,P\31_{\Gamma\backslash Y_J}.
$$
Using formula~(\ref{epro}) for the set $J$ instead of $F$, we also
see that $P_J\31_{\Gamma\backslash Y_J}$ is the constant function
$\zeta_{S_J,\Gamma}(\beta)^{-1}$. Using again that $PP_J=P$, we
therefore obtain
\begin{equation} \label{epro2}
Pf_J=\zeta_{S_J,\Gamma}(\beta)^{-1}\int_{\Gamma\backslash\pgl}
fd\bar\nu_\infty.
\end{equation}
Since the space $H$ contains nonzero constant functions, this in
particular implies that
$$
\int f_Jd\tilde\nu=\zeta_{S_J,\Gamma}(\beta)^{-1}
\int_{\Gamma\backslash\pgl}fd\bar\nu_\infty,
$$
so that $\int f_Jd\tilde\nu$ is the same for every $\tilde\mu$.

\smallskip

To extend the result to all $\G(\Zhat)$-invariant functions, fix a
finite nonempty set $J$ of prime numbers, and consider a right
$\prod_{p\in J}\G(\Z_p)$-invariant bounded Borel function $f$ on
$$
\Gamma\backslash\left(\pgl\times\prod_{p\in J}\mtwo(\Z_p)\right).
$$
We may consider $f$ as a function on
$\Gamma\backslash(\pgl\times\mtwo(\Zhat))$. Then $f$ is right
$\G(\Zhat)$-invariant, and the space spanned by such functions for all
$J$'s is dense in the space of square integrable $\G(\Zhat)$-invariant
functions. Applying again formula~(\ref{epro}) for the projection
$P_J$ (for $J$ in place of~$F$), we see that $P_Jf$ is again a
function whose value at $(\tau,m)\in\pgl\times\mtwo(\Zhat)$ depends
only on~$\tau$ and~$m_p$ with $p\in J$. The formula also shows that
$P_J$ commutes with the action of $\G(\Zhat)$, so $P_Jf$ is
$\G(\Zhat)$-invariant. Since $\G(\Z_p)$ acts transitively on itself,
this shows that the value of $P_Jf$ at $(\tau,m)$ with $m_p\in
\G(\Z_p)$ for $p\in J$ depends only on $\tau$. In other words, on
the space $\Gamma\backslash Y_J$ introduced above, the function
$P_Jf$ is a bounded Borel function of the form $\tilde f_J$ for some
function $\tilde f$ on $\Gamma\backslash\pgl$. An important point is
that $\tilde f$ depends on $f$ but not on $\tilde\mu$. By
Lemma~\ref{project}(1) and the polarization identity we have
$$
\int
P_Jfd\tilde\nu=\zeta_{S_J,\Gamma}(\beta)\int_{\Gamma\backslash
Y_J}P_Jfd\tilde\nu=\zeta_{S_J,\Gamma}(\beta)
\int_{\Gamma\backslash Y_J} \tilde
f_Jd\tilde\nu=\int_{\Gamma\backslash\pgl}\tilde fd\bar\nu_\infty.
$$
Since $\int fd\tilde\nu=\int P_Jfd\tilde\nu$, we see again that
$\int fd\tilde\nu$ is the same for any $\tilde\mu$. It therefore
follows that $\int fd\tilde\nu=\int fd\bar\nu_\beta$ for any
bounded Borel $\G(\Zhat)$-invariant function on
$\Gamma\backslash(\pgl\times\mtwo(\Zhat))$. Since~$\tilde\nu$ is
$\G(\Zhat)$-invariant by assumption, we have
$\tilde\nu=\bar\nu_\beta$ and hence $\tilde\mu=\bar\mu_\beta$.

\smallskip

To prove ergodicity assume $Z_0$ is a left $\glq$-invariant and
right $\glr$-invariant $\bar\mu_\beta$-measur\-able subset of
$\pgl\times\ma$ of positive measure. Since
$\glq(\pgl\times\mr)=\pgl\times\ma$, it follows that the set
$Z_0\cap(\pgl\times\mr)$ has positive measure. Hence
$\lambda=\bar\mu_\beta(\Gamma\backslash(Z_0\cap(\pgl\times\mr)))>0$.
It follows that the measure $\tilde\mu$ defined by
$$
\tilde\mu(Z)=\lambda^{-1}\bar\mu_\beta(Z_0\cap Z)
$$
is right $\glr$-invariant and satisfies the conditions in
Theorem~\ref{pgl}. Hence $\tilde\mu=\bar\mu_\beta$, and consequently
the complement of $Z_0$ has $\bar\mu_\beta$-measure zero. \ep

We aim to prove that the action of $\glq$ on
$(\pgl\times\ma,\bar\mu_\beta)$ is ergodic. The next step is to
consider the action on $\ma$ alone.

\begin{lemma} \label{finite}
The action of $\glq$ on $(\ma,\mu_{\beta,f})$ is ergodic.
\end{lemma}

\bp The proof is similar to that of the previous lemma, but
requires a much simpler result than Lemma~\ref{clozel}.

Consider the space $L^2(\mtwo(\Zhat),d\mu_{\beta,f})$ and the
subspace $H$ of $\mzp$-invariant functions. It suffices to show
that $H$ consists of constant functions. Denote by $P$ the
orthogonal projection onto~$H$.

For a finite set $F$ of prime numbers denote by $P_F$ the
projection onto the space of $S_F$-invariant functions. Put also
$$
Y_F=\prod_{p\in F}\G(\Z_p)\times\prod_{q\notin F}\mtwo(\Z_q).
$$
Then similarly to~(\ref{epro}) for any $\Gamma$-invariant function
$f\in L^2(\mtwo(\Zhat),d\mu_{\beta,f})$ we have
\begin{equation} \label{epro4}
P_Ff|_{S_Fm}=\zeta_{S_F,\Gamma}(\beta)^{-1}\sum_{s\in\Gamma
\backslash S_F/\Gamma}\det(s)^{-\beta}R_\Gamma(s)(T_sf)(m)\ \
\text{for}\ \ m\in Y_F.
\end{equation}
This can be either proved similarly to Lemma~\ref{project}(2) or
deduced from that lemma by identifying the space of
$\Gamma$-invariant functions with the subspace of
$L^2(\Gamma\backslash(\pgl\times\mtwo(\Zhat)),d\bar\nu_\beta)$ of
functions depending only on the second coordinate.

For a finite set $J$ of primes disjoint from $F$, and a left
$\Gamma$-invariant function $f$ on $\prod_{p\in J}\G(\Z_p)$ define a
function $f_J$ by
$$
f_J(m)=\begin{cases}
f((m_p)_{p\in J}), & \text{if}\ m_p\in\G(\Z_p)\ \text{for}\ p\in J,\\
0, &\text{otherwise.}\end{cases}
$$
Since $f$ is $\Gamma$-invariant and $\Gamma$ is dense in
$\prod_{p\in J}\sltwo(\Z_p)$, $f$ is invariant with respect to
multiplication on the left by elements of the latter group. In
other words, the value of $f$ at $m$ depends only on
$\det(m)\in\prod_{p\in J}\Z_p^*$. Therefore functions of the form
$f(m)=\chi(\det(m))$, where $\chi$ is a character of the compact
abelian group $\prod_{p\in J}\Z_p^*$, span a dense subspace of
$\Gamma$-invariant functions on $\prod_{p\in J}\G(\Z_p)$. But if
$f=\chi\circ\det$, we have
$$
(T_sf)(m)=\chi(\det(m))\chi(\det(s))
$$
for $s\in S_F$ and $m\in\prod_{p\in J}\G(\Z_p)$. Applying now
(\ref{epro4}) to the function $f_J$ and using a calculation similar
to~(\ref{ezeta}) and~(\ref{ezeta1}), we get
\begin{eqnarray*}
P_Ff_J|_{S_Fm}&=&\chi(\det((m_p)_{p\in
J}))\zeta_{S_F,\Gamma}(\beta)^{-1} \sum_{s\in\Gamma \backslash
S_F/\Gamma}\det(s)^{-\beta}R_\Gamma(s)\chi(\det(s))\\
&=&\chi(\det((m_p)_{p\in J}))\prod_{p\in
F}\frac{(1-p^{-\beta})(1-p^{-\beta+1})}
{(1-\chi(p)p^{-\beta})(1-\chi(p)p^{-\beta+1})}.
\end{eqnarray*}
If the character $\chi$ is nontrivial, by choosing $F$ large enough
the product above can be made arbitrarily small by elementary
properties of Dirichlet series (this was used already for the
classification of KMS-states of the Bost-Connes system
in~\cite{bos-con}, see also~\cite{nes}). Since $PP_F=P$, we conclude
that $Pf_J=0$. On the other hand, if $\chi$ is trivial then
$f_J=\31_{Y_J}$. Then applying~(\ref{epro4}) with $J$ in place of
$F$ we get $P_Jf_J=\zeta_{S_J,\Gamma}(\beta)^{-1}$. In either case
we see that $Pf_J$ is constant.

\smallskip

Let now $f$ be a function on $\prod_{p\in J}\G(\Z_p)$ which is no
longer left $\Gamma$-invariant. Since $\Gamma$ is dense
in~$\sltwo(\Zhat)$, any function in $H$ is $\sltwo(\Zhat)$-invariant.
Hence to compute $Pf_J$ we can first apply to $f_J$ the projection
$Q$ onto the subspace of $\sltwo(\Zhat)$-invariant functions. But $Q$
is given by averaging over $\sltwo(\Zhat)$-orbits. We then see that
$Qf_J=\tilde f_J$, where
$$
\tilde f(m)=\int_{\prod_{p\in J}\sltwo(\Z_p)}f(gm)dg.
$$
Hence $Pf_J=PQf_J=P\tilde f_J$ is again a constant function.

\smallskip

To extend the result to all functions on $\mr$, for each $s\in\mzp$
we introduce an operator~$V_s$ on the space
$L^2(\mtwo(\Zhat),d\mu_{\beta,f})$ by letting $(V_sh)(m)=h(sm)$. Then
$V_sP=P$. Using the scaling condition we see that
$\det(s)^{-\beta/2}V_s$ is a coisometry with initial space
$L^2(s\mtwo(\Zhat),d\mu_{\beta,f})$. It follows that the adjoint
operator is given by
$$
(V^*_sh)(y)=\begin{cases}\det(s)^{\beta}h(s^{-1}y),&\text{if}\ y\in
s\mtwo(\Zhat),\\ 0,&\text{otherwise}.\end{cases}
$$
In particular, we see that if $s\in S_J$ for some finite set $J$
then both operators $V_s$ and $V^*_s$ preserve the space of
functions $f$ such that $f(m)$ depends only on $m_p$ with $p\in J$.
But then if $f$ is such a function with support on $Y_J$, the
function $V^*_sf$ has support on $sY_J$. Since $V_sP=P$, we have
$PV_s^*=P$ and thus $PV_s^*f=Pf$ is a constant. We thus see that the
image of a dense space of functions consists of constant functions.
\ep

The following simple trick will allow us to combine the two previous
lemmas. It expounds a remark in~\cite{nes}.

\begin{proposition} \label{comm}
Assume we have mutually commuting actions of locally compact
second countable groups $G_1$, $G_2$ and $G_3$ on a Lebesgue space
$(X,\mu)$. Suppose that \enu{i} the actions of $G_1$ on
$(X/G_2,\mu)$ and $(X/G_3,\mu)$ are ergodic; \enu{ii} $G_2$ is
connected and $G_3$ is compact totally disconnected.

Then the action of $G_1$ on $(X,\mu)$ is ergodic.
\end{proposition}

Here by quotient spaces we mean quotients in measure theoretic
sense. So by definition
$$
L^\infty(X/G_i,\mu)=L^\infty(X,\mu)^{G_i}.
$$

\bp[Proof of Proposition~\ref{comm}] By assumption the action of
$G_1\times G_3$ on~$X$ is ergodic. In other words, the action of
$G_3$ on~$X/G_1$ is ergodic. Since $G_3$ is compact, we can then
identify $X/G_1$ with a homogeneous space of~$G_3$, say $G_3/H$,
where $H$ is a closed subgroup of~$G_3$, see e.g.
\cite[Section~2.1]{zim}. Since $G_1\times G_2$ acts ergodically on
$X$, we have an ergodic action of~$G_2$ on~$X/G_1=G_3/H$. Since this
action commutes with the action of $G_3$ on $G_3/H$ by left
translations, it is given by right translations, that is, by a
measurable homomorphism $G_2\to N(H)/H$, where $N(H)$ is the
normalizer of~$H$ in~$G_3$. Such a homomorphism is automatically
continuous (see e.g.~\cite[Theorem~B.3]{zim}), and since $G_2$ is
connected and $N(H)/H$ is totally disconnected, the homomorphism
must be trivial. But since the action of~$G_2$ is ergodic this means
that $H=G_3$, so that $X/G_1$ is a single point. Thus the action
of~$G_1$ is ergodic. \ep

\begin{corollary}\label{glqergodicpglma}
The left action of $\glq$ on $(\pgl\times\ma, \bar\mu_\beta)$ is
ergodic.
\end{corollary}

\bp The proof  is a straightforward application of
Proposition~\ref{comm} with $G_1=\glq$, $G_2=\pgl$ and
$G_3=\G(\Zhat)$, so that $G_1$ acts on $\pgl\times\ma$ by
multiplication on the left and $G_2$ and $G_3$ act by multiplication
on the right on the corresponding factor. That the actions of~$G_1$
on the quotient spaces are ergodic is given by Lemma~\ref{symm} and
Lemma~\ref{finite}. \ep

\bp[Proof of Theorem~\ref{pgl}] We follow an argument similar to
that of ~\cite[Theorem~25]{bos-con}. Note first that $\bar\mu_\beta$
is right $\G(\Zhat)$-invariant and satisfies  the conditions in
Theorem~\ref{pgl} by construction. Denote by $K_\beta$ the affine
set of measures on $\pgl\times\ma$ satisfying the conditions in
Theorem~\ref{pgl}. Let $C_\beta$ be a cone with base $K_\beta$.
Denote by $v_0$ its vertex. The cone $C_\beta$ has the structure of
a Choquet simplex. Namely, similarly to Proposition~\ref{cm1} it can
be identified with the set of KMS$_\beta$-states on~$B^\sim$, where
$$
B=\red{\fibb{\glq}{\Gamma}{(\pgl\times\mtwo(\Zhat))}},
$$
$B^\sim$ is obtained from~$B$ by adjoining a unit, and $v_0$
corresponds to the state on $B^\sim$ with kernel $B$. Denote by
$\bar\varphi$ the state corresponding to~$\bar\mu_\beta$. Then by
Remark~\ref{vonneumann} the algebra $\pi_{\bar\varphi}(B^\sim)''$ is
a reduction of the von Neumann algebra of the orbit equivalence
relation defined by the action of $\glq$ on
$(\pgl\times\ma,\bar\mu_\beta)$. By Corollary~\ref{glqergodicpglma}
this von Neumann algebra is a factor. Hence
$\pi_{\bar\varphi}(B^\sim)''$ is also a factor, and therefore
$\bar\mu_\beta$ is an extremal point of $C_\beta$. The group
$\G(\Zhat)$ acts on~$C_\beta$, and by virtue of Lemma~\ref{symm} the
segment $[\bar\mu_\beta,v_0]$ is the set of $\G(\Zhat)$-invariant
points. Suppose now  $v \in C_\beta$ is an extremal point. Then $w =
\int_{\glr} gv\, dg\in C_\beta$ where  each $gv$ is also an extremal
point of~$C_\beta$. But because of its $\glr$ invariance, $w$ lies
on $[\bar\mu_\beta,v_0]$ and hence is also a convex combination of
the extremal points $\bar\mu_\beta$ and $v_0$. Since $w$ is the
barycenter of a unique probability measure on the set of extremal
points, we conclude that either $v=\bar\mu_\beta$ or $v=v_0$. Thus
$C_\beta=[\bar\mu_\beta,v_0]$ and $K_\beta=\{\bar\mu_\beta\}$. This
completes the proof of Theorem~\ref{pgl}. \ep

\begin{remark} \label{rsingular}
We have classified KMS$_\beta$-states of the Connes-Marcolli system
for $\beta\ne0,1$. Let us now briefly discuss the cases $\beta=0,1$.
\enu{i}If $\beta=0$ then by Lemma~\ref{singular} and the
considerations following Corollary~\ref{singular2} one can conclude
that there are no nonzero finite traces on
$I=\red{\fibb{\glq}{\Gamma}{(\HH\times\mtwo(\Zhat)^\times)}}$.
Therefore the only KMS$_0$-states, that is, $\sigma$-invariant
traces, are those coming from $A/I=\red{\fib \glq \Gamma \HH}$.
There is a canonical trace defined by the $\glq$-invariant measure
$\mu_\infty$ on $\HH$. Notice that though the action of $\glq$ on
$\HH$ is not free and so Proposition~\ref{Scaling1} is not
immediately applicable, the action of $\glq/\Q^*$ is free in the
measure theoretic sense, and this is enough to check the trace
property. This is probably the unique such trace. \enu{ii} If
$\beta=1$ then, as we know, KMS$_1$-states still correspond to
measures satisfying the scaling condition. By the first part of the
proof of Corollary~\ref{singular2} and our considerations following
that corollary, the set of points $(\tau,m)\in\HH\times\ma$ with
$m_p\ne0$, $\det(m_p)=0$ for every $p$, is a subset of full measure.
Such measures indeed exist. Let $\mu'_f$ be the Haar measure on the
locally compact group $\af^2$ normalized such that
$\mu'_f(\Zhat^2)=1$. We may consider $\mu'_f$ as a measure on $\ma$ by
identifying $\af^2$ with the set of matrices with zero first column.
Then $\mu'=2\mu_\infty\times\mu'_f$ is a measure with the required
properties. Using the action of $\G(\Zhat)$ by multiplication on the
right we can then construct infinitely many such measures (notice
that the stabilizer of $\mu'$ in $\G(\Zhat)$ is the group of upper
triangular matrices). We conjecture that this way one gets all
extremal KMS$_1$-states.
\end{remark}

\begin{remark}
Let $1<\beta\le2$, and denote by $\varphi_\beta$ the unique
KMS$_\beta$-state on the Connes-Marcolli C$^*$-algebra $A$. It is
easy to describe the flow of weights of the factor
$\pi_{\varphi_\beta}(A)''$. Let us first consider the algebra
$B=\red{\fibb{\glq}{\Gamma}{(\pgl\times\mtwo(\Zhat))}}$ and the state
$\bar\varphi_\beta$ on~$B$ corresponding to~$\bar\mu_\beta$, and
describe the flow of weights of $\pi_{\bar\varphi_\beta}(B)''$. By
Remark~\ref{vonneumann}, equivalently we want to describe the flow
of weights of the orbit equivalence relation defined by the ergodic
action of $\glq$ on $(\pgl\times\ma,\bar\mu_\beta)$.

The group $\R^*_+$ acts on the measure space
$(\R^*_+\times\pgl\times\ma,
\lambda\times\bar\mu_\infty\times\mu_{\beta,f})$, where $\lambda$
is a measure in the Lebesgue measure class, by
$$
t(s,h,\rho)=(t^{-1/\beta}s,h,\rho).
$$
The flow of weights is induced by this action on the quotient of
the space by the action of $\glq$ defined by
$$
g(s,h,\rho)=(\det(g)s,gh,g\rho).
$$
We have an isomorphism $\glp(\R)/ \{\pm1\}\to\R^*_+\times\pgl$,
$g\mapsto(\det(g),\bar g)$, where $\bar g$ denotes the class of $g$
in $\pgl$. So instead of the space $\R^*_+\times\pgl\times\ma$ we
may consider $(\glp(\R)/ \{\pm1\})\times\ma$. We may further replace
$\glp(\R)$ by $\G(\R)$, but instead of the action of $\glq$ we then
have to consider the action of $\G(\Q)$. Finally, replace $\G(\R)$
by $\mtwo(\R)$, and so instead of $\G(\R)\times\ma$ consider
$\mtwo(\A)$, where $\A=\R\times\A_f$ is the full adele space. To
summarize, $\R^*_+$ acts on $\mtwo(\A)=\mtwo(\R)\times\ma$ by
$t(m,\rho)=(t^{-1/2\beta}m,\rho)$, and the flow of weights of the
factor $\pi_{\bar\varphi_\beta}(B)''$ is induced by this action on
the quotient of the measure space
$(\mtwo(\R)\times\ma,\lambda_\infty\times\mu_{\beta,f})$, where
$\lambda_\infty$ is the usual Lebesgue measure on
$\mtwo(\R)\cong\R^4$, by the action of $\G(\Q)\times\{\pm1\}$
defined by $(g,s)(m,\rho)=(gms,g\rho)$.

Denote the measure $\lambda_\infty\times\mu_{\beta,f}$ on
$\mtwo(\A)$ by $\lambda_\beta$. Note that $\lambda_2$ is a Haar
measure on the additive group $\mtwo(\A)$.

Similarly, by identifying $\R^*_+\times\HH$ with $\glp(\R)/\so$ we
conclude that the flow of weights of the factor
$\pi_{\varphi_\beta}(A)''$ is defined on the quotient of the measure
space $(\mtwo(\A),\lambda_\beta)$ by the action of $\G(\Q)\times\so$
defined by $(g,s)(m,\rho)=(gms,g\rho)$ for $(g,s)\in\G(\Q)\times\so$
and $(m,\rho)\in\mtwo(\A)=\mtwo(\R)\times\ma$.

It seems natural to conjecture that the action of $\G(\Q)$ on
$(\mtwo(\A),\lambda_\beta)$ is ergodic, so the flows of weights of
the factors $\pi_{\varphi_\beta}(A)''$ and
$\pi_{\bar\varphi_\beta}(B)''$ are trivial, and thus the factors are
of type~III$_1$. The analogous property in the one-dimensional case
indeed holds~\cite{bos-con,nes}. Note that so far we have only shown
that the action of $\G(\Q)\times\R^*$ is ergodic, which is
equivalent to ergodicity of the action of $\glq$ on
$(\pgl\times\ma,\bar\mu_\beta)$. Note also that similarly to the
one-dimensional case~\cite{nes}, by virtue of Lemma~\ref{finite} and
Proposition~\ref{comm}, to prove the conjecture it would be enough
to show that the action of~$\glq$ on $\glp(\R)\times(\ma/\G(\Zhat))$
is ergodic, or equivalently, the action of $\G(\Q)$ on
$(\mtwo(\A)/\G(\Zhat),\lambda_\beta)$ is ergodic.  Recall that in the
one-dimensional case the corresponding ergodicity result for the
action of $\Q^*$ on $\A/\Zhat^*$ was established in~\cite{bla}
and~\cite{boc-zah}.
\end{remark}

\begin{remark} \label{rgln}
We believe that the results of Sections~\ref{scm} and~\ref{scm2}
are valid for $\gln$ for any $n\ge2$. More precisely, consider the
algebra $\red{\fibb{\glnq}{\slnz}{(\pgln\times\mn(\Zhat)})}$. Define
a dynamics by the homomorphism $\glnq\ni g\mapsto\det(g)$. Then
\enu{i} for $\beta\in(-\infty,0)\cup(0,1)\cup\dots\cup(n-2,n-1)$
there are no KMS$_\beta$-states; \enu{ii} for $\beta\in(n-1,n]$
there exists a unique KMS$_\beta$-state; \enu{iii}\, for $\beta>n$
there is a one-to-one correspondence between KMS$_\beta$-states
and probability measures on
$\slnz\backslash(\pgln\times\gln(\Zhat))$; \enu{iv}\, for
$\beta=0,1,\dots,n-1$ there is a KMS$_\beta$-state defined by the
Haar measure on $\af^{\beta n}$, when we identify the latter group
with the set of matrices in $\mn(\af)$ with zero first $n-\beta$
columns.

The key step for this generalization would be an analogue of
Lemma~\ref{singular}.
\end{remark}

\end{document}